\documentclass[11pt, reqno]{amsart}
\usepackage{amsfonts}
\usepackage{amsmath}
\usepackage{color}
\usepackage{latexsym}
\usepackage{amssymb}
\usepackage{mathrsfs}
\usepackage{amsmath}
\usepackage{fancybox,color}
\usepackage{enumerate}
\usepackage[active]{srcltx}
\usepackage[colorlinks]{hyperref}
\usepackage{hyperref}
\usepackage{hypernat}
\usepackage{eurosym}
\usepackage{url}
\definecolor{blau}{rgb}{0.1,0.0,0.9}
\definecolor{purple}{rgb}{0.4,0.0,0.9} 
\definecolor{gruen}{cmyk}{1.0,0.2,0.7,0.07}
\definecolor{mag}{cmyk}{0.0,0.9,0.3,0.0}

\newcommand{\kom}[1]{}
\newcounter{komcounter}
\numberwithin{komcounter}{section}

 \def\1{\raisebox{2pt}{\rm{$\chi$}}}
\def\a{{\bf a}}

\newtheorem{theorem}{Theorem}[section]
\newtheorem{corollary}[theorem]{Corollary}
\newtheorem{lemma}[theorem]{Lemma}
\newtheorem{proposition}[theorem]{Proposition}
\newtheorem{definition}[theorem]{Definition}
\newtheorem{remark}[theorem]{Remark}

\newcommand{\R}{{\mathbb R}}
\newcommand{\re}{{\mathbb R}}

\newcommand{\N}{{\mathbb N}}

\newcommand{\B}{{\mathcal B^{q,d,p}}}

 \newcommand{\eps}{{\varepsilon}}
 \def\1{\raisebox{2pt}{\rm{$\chi$}}}
 
\newcommand{\abs}[1]{\left|#1\right|}
\newcommand{\norm}[1]{\left|\left|#1\right|\right|}
\newcommand{\Rn}{\mathbb{R}^n}
\newcommand{\ren}{\mathbb{R}^n}
\newcommand{\osc}{\operatorname{osc}}

\def\vint_#1{\mathchoice%
          {\mathop{\kern 0.2em\vrule width 0.6em height 0.69678ex depth -0.58065ex
                  \kern -0.8em \intop}\nolimits_{\kern -0.4em#1}}%
          {\mathop{\kern 0.1em\vrule width 0.5em height 0.69678ex depth -0.60387ex
                  \kern -0.6em \intop}\nolimits_{#1}}%
          {\mathop{\kern 0.1em\vrule width 0.5em height 0.69678ex
              depth -0.60387ex
                  \kern -0.6em \intop}\nolimits_{#1}}%
          {\mathop{\kern 0.1em\vrule width 0.5em height 0.69678ex depth -0.60387ex
                  \kern -0.6em \intop}\nolimits_{#1}}}
\def\vintslides_#1{\mathchoice%
          {\mathop{\kern 0.1em\vrule width 0.5em height 0.697ex depth -0.581ex
                  \kern -0.6em \intop}\nolimits_{\kern -0.4em#1}}%
          {\mathop{\kern 0.1em\vrule width 0.3em height 0.697ex depth -0.604ex
                  \kern -0.4em \intop}\nolimits_{#1}}%
          {\mathop{\kern 0.1em\vrule width 0.3em height 0.697ex depth -0.604ex
                  \kern -0.4em \intop}\nolimits_{#1}}%
          {\mathop{\kern 0.1em\vrule width 0.3em height 0.697ex depth -0.604ex
                  \kern -0.4em \intop}\nolimits_{#1}}}

\newcommand{\kint}{\vint}

\newcommand{\aveint}[2]{\mathchoice%
          {\mathop{\kern 0.2em\vrule width 0.6em height 0.69678ex depth -0.58065ex
                  \kern -0.8em \intop}\nolimits_{\kern -0.45em#1}^{#2}}%
          {\mathop{\kern 0.1em\vrule width 0.5em height 0.69678ex depth -0.60387ex
                  \kern -0.6em \intop}\nolimits_{#1}^{#2}}%
          {\mathop{\kern 0.1em\vrule width 0.5em height 0.69678ex depth -0.60387ex
                  \kern -0.6em \intop}\nolimits_{#1}^{#2}}%
          {\mathop{\kern 0.1em\vrule width 0.5em height 0.69678ex depth -0.60387ex
                  \kern -0.6em \intop}\nolimits_{#1}^{#2}}}
\newcommand{\ud}{\, d}
\newcommand{\half}{{\frac{1}{2}}}
\newcommand{\parts}[2]{\frac{\partial {#1}}{\partial {#2}}}
\newcommand{\ol}{\overline}
\newcommand{\Om}{\Omega}

\newcommand{\dist}{\operatorname{dist}}

\newcommand{\om}{\omega}

\newcommand{\spt}{\operatorname{spt}}
\newcommand{\vp}{\varphi}

\newcommand{\sgn}{\operatorname{sgn}}

\renewcommand{\a}{\alpha}

\newcommand{\ve}{\varepsilon}

\def\qed{\,\unskip\kern 6pt \penalty 500
\raise -2pt\hbox{\vrule \vbox to8pt{\hrule width 6pt
\vfill\hrule}\vrule}\par}

\numberwithin{equation}{section}
\begin{document}

\title[]{Equivalence between radial solutions of different parabolic gradient-diffusion equations and applications}

\author[]{Mikko Parviainen, Jyv\"askyl\"a\\ Juan Luis V\'azquez,  Madrid}
\date{\today}

\maketitle

\begin{abstract}
We consider a general form of a parabolic  equation  that generalizes both  the standard parabolic  $p$-Laplace equation and the normalized version that has been proposed in stochastic game theory.
We establish an equivalence between this equation and the standard $p$-parabolic
equation posed in a fictitious space dimension, valid for radially symmetric solutions.
This allows us to find suitable explicit solutions for example of Barenblatt type, and as a consequence we settle  the exact asymptotic behaviour of the Cauchy problem even for nonradial data.  We also establish the asymptotic behaviour in a bounded domain. Moreover, we use the explicit solutions to establish the parabolic Harnack's inequality.
\end{abstract}

\footnotesize

\normalsize

\section{Introduction}
\label{sec:intro}

In this paper we consider the  following general version of the nonlinear parabolic equation
\begin{equation}\label{eq:equation}
\partial_t u= c \abs{Du}^{\kappa}\operatorname{div}(\abs{Du}^{p-2} Du),
\end{equation}
with real parameters  $p>1$ and $\kappa>1-p$,  so that the homogeneity of the right-hand side operator is always larger than 0. The constant $c>0$ does not play any role for the given equation since it can be eliminated by a time scaling, but it could be useful to play with when the exponents vary, mainly when $p\to \infty$. Formally \eqref{eq:equation} is a nonlinear parabolic equation, possibly degenerate or singular; moreover, for $\kappa\ne 0$ the right-hand side is not a divergence-form operator. The equation gives the usual $p$-parabolic equation and also the  normalized $p$-parabolic equation, both will be briefly described below as a way of motivating the general problem. Note that for $p=2$ we get the equation $\partial_t u=\abs{Du}^{\kappa}\Delta u$.

We will discuss the initial-value problem in a bounded domain  $\Omega\subset \ren$, or  in the whole space $\ren$, $n\ge 2$. Our main technical contribution is the reduction of the general equation to the standard $p$-parabolic equation by an equation transformation. The equivalence applies only to radially symmetric solutions but this will be enough to find suitable special solutions of Barenblatt type. Moreover, this helps us to settle the long-time behaviour of all solutions  with continuous, bounded and compactly supported initial data, as well as to establish the parabolic Harnack's inequality.

\medskip

\noindent {\sc Standard parabolic $p$-Laplace equation.} The most popular model of evolution equation of nonlinear diffusion type with gradient-dependent diffusivity is the so-called {\sl parabolic $p$-Laplace equation} (also known as $p$-parabolic equation for short):
\begin{equation}
\nonumber
\partial_t u= \Delta_p u:= \operatorname{div}(\abs{Du}^{p-2} Du)\,.
\end{equation}
It has been extensively studied for all  values of the parameter $p\in (1,\infty)$. We use the notation $Du(x,t)$ to denote the spatial gradient of functions $u(x,t)$, $x\in \ren$, $n\ge 1$, and $ t\in \re$. Note that for $p=2$ we recover the classical heat equation. For $p>2$ the diffusivity $\abs{Du}^{p-2}$ vanishes for $Du=0$ and therefore this is called degenerate or slow diffusion, while for $1<p<2$ we have $\abs{Du}^{p-2}\to \infty$ as $\abs{Du}\to 0$ and this is called singular or fast diffusion.

The parabolic $p$-Laplace equation has a large literature covering existence and uniqueness of weak solutions  for suitable initial and boundary data, for problems posed in the whole space or in bounded domains. Also regularity, asymptotic behaviour, free boundaries and other issues have been studied, see for example the monograph \cite{dibenedetto93}, and more recently \cite{vazquez06, urbano08}.
These works use the theory of weak solutions. Since the equation (\ref{eq:equation}) is in non-divergence form except in this particular  special case, the solutions in this paper are understood in the viscosity sense, see \cite{crandallil92, giga06,ohnumas97}. However, in the case of $p$-parabolic equation the notions of weak  and viscosity solution are equivalent for all $p\in (1,\infty)$ \cite{juutinenlm01}.

\medskip

\noindent {\sc The limit cases.} Also the limit cases $p=1$  and $p=\infty$  have attracted attention, and have posed problems that help  motivate our work.
The case $p=1$ is a model for the so called {\sl total variation flow}
\begin{equation}
\nonumber
\partial_t u= \Delta_1 u:= \operatorname{div}\left(\frac{ Du}{\abs{Du}} \right)\,
\end{equation}
and this equation appears for example in image processing, see \cite{andreybcm01}. However, there is another model that involves the 1-Laplace operator, namely, the mean curvature flow
\begin{equation}
\nonumber
\partial_t u= \abs{Du}\,\Delta_1 u:=\abs{Du}\,\operatorname{div}\left(\frac{ Du}{\abs{Du}}\right)\,,
\end{equation}
which is very important in differential geometry \cite{evanss91, chengg91}, and in a number of applications, like crystal growth. Here we find a first case of occurrence of the non-divergence factor $\abs{Du}^{\kappa}$, with $\kappa=1$.

The limit as $p\to\infty$ posed another problem  \cite{juutinenk06} and opened up another perspective.  For a smooth $p$-harmonic function with nonvanishing gradient, a short calculation shows that
$$
\Delta_p u=|Du|^{p-2} \, \left(\Delta u + (p-2) \sum_{i,j=1}^n \frac{\partial_{i}u}{|Du|}\frac{\partial_{j}u}{|Du|}\partial_{ij}u \right)=0.
$$
This formula suggests that the factor $|Du|^{p-2}$ could be disregarded (in this particular instance), and attention could be concentrated on the operator  between parentheses. This can be made rigorous, since in the viscosity formulation of the $p$-Laplace equation, we can disregard the test functions with vanishing gradient \cite{juutinenlm01}.
Further,  formally dividing by $p$ and passing to the limit $p\to\infty$ produces the  normalized or game theoretic version of the infinity Laplace operator
\begin{equation}
\nonumber
\Delta_{\infty}^{N} u=\sum_{i,j=1}^n \frac{\partial_{i}u}{|Du|}\frac{\partial_{j}u}{|Du|}\partial_{ij}u\,.
\end{equation}
Next we observe that eliminating factors of the form $|Du|^\gamma$ is not possible in the evolution problem because of the presence of the left-hand side $\partial_t u$. Therefore, the standard version of the $p$-Laplace evolution equation can be written as
\begin{equation}
\nonumber
\partial_t u= \abs{Du}^{p-2} \left(\Delta u + (p-2)\Delta_{\infty}^{N}u \right)\,,
\end{equation}
and it is not equivalent to the modified equation
\begin{equation}
\label{eq:normalized}
\partial_t u=  \Delta u + (p-2)\Delta_{\infty}^{N}u=:\Delta_p^N u\,.
\end{equation}
The new $p$-Laplace operator $\Delta_p^N u$, called the {\sl normalized or game theoretic $p$-Laplacian}, can be seen as an interpolation between the standard Laplace and $\Delta_{\infty}^{N}$, and has been proposed in stochastic game theory as explained below.

\medskip

\noindent {\sc Viscosity solutions.} The works of Crandall, Evans, Giga, Ishii, Lions, Souganidis and others established the  theory of  viscosity solutions and their connection to the stochastic differential games and control theory in the early 80s.
Recently, a connection between the theory of stochastic tug-of-war games and normalized $p$-Laplace type equations
$$
\Delta_p^N u=f
$$
 has been investigated. In the elliptic case $1<p\le \infty$, this connection was discovered  in the
works of Peres, Schramm, Sheffield and Wilson \cite{peress08, peresssw09}. For $p=1$, see \cite{buckdahncq01,kohns06}.  In the parabolic case, it was shown in \cite{manfredipr10} that solutions to (\ref{eq:normalized}) can be obtained as limits of values of tug-of-war games with noise when the parameter that controls the length of steps goes to zero.
Solutions to equations of type (\ref{eq:normalized}) remain solutions when multiplied by a constant, which can be a useful attribute in applications to mathematical image processing \cite{does11,elmoataztt15}: then the brightness of the original image does not affect the evolution itself. Heuristically speaking, the value of $p$ in (\ref{eq:normalized})  controls the strength of the diffusion to the direction of sharp changes (gradient direction) compared to the diffusion to other directions when processing images using this evolution.

Recently,  equation (\ref{eq:normalized}) has been studied for example by Banerjee and Garofalo \cite{banerjeeg13, banerjeeg15}, Juutinen \cite{Ju14}, Jin and Silvestre \cite{jins17}, Attouchi and Parviainen \cite{attouchip}, Ubostad \cite{ubostad}, Berti and Magnanini \cite{bertim} as well as A.~Bj\"orn, J.~Bj\"orn and Parviainen \cite{bjornbp}. The general equation (\ref{eq:equation}) can be badly singular and thus requires a modification of the definition of viscosity solutions as in \cite{ohnumas97}, see also \cite{ishiis95} and \cite{demengel11}. $C^{1,\a}$-regularity for (\ref{eq:equation}) was established by Imbert, Jin and Silvestre in \cite{imbertjs}.

\medskip

\noindent {\sc Outline of the paper.} The above considerations motivate the study of the more general form \eqref{eq:equation}.  In Section \ref{sec:visc} we introduce a suitable concept of viscosity solution taken from Ohnuma and Sato \cite{ohnumas97}, and recall existence, uniqueness and comparison results.

In Section \ref{sec:radial} we present a transformation that works for radially symmetric solutions and allows to reduce the general case of equation \eqref{eq:equation}  to a version of the known $p$-Laplace theory or its radial counterpart. Heuristically speaking, we show that the solutions to the original problem can be interpreted as solutions to the divergence form $p$-parabolic equation, but in an fictitious space dimension $d$ given by
\begin{equation}
d-1=(n-1)\,\frac{\kappa+p-1}{p-1}\,.
\end{equation}
This is inspired by the works \cite{portilheirov12,portilheirov13} about  evolution flows driven by the infinity Laplacian, where $p=\infty$ and $d=1$. However, since $d$ is not necessarily integer, the rigorous connection  for $p<\infty$  requires looking at a weighted 1-dimensional divergence form equation, Definition \ref{def:1D-weak}, for which the connection is established in Theorem \ref{thm:equiv}.

The equivalence result then allows us to derive important explicit radial solutions to  \eqref{eq:equation} in the  examples  of Subsection \ref{ss.examples}: the exponential solutions,  the Barenblatt type solutions; and later in Section  \ref{sec:bounded-dom} we introduce the friendly giant type solutions.

As an application, using Aleksandrov's reflection principle and the scaling properties of the Barenblatt solution, we show in Theorem \ref{thm:Cauchy-asymp} that the viscosity solution to  equation (\ref{eq:equation})  behaves asymptotically like a Barenblatt solution of the  preceding $p$-Laplacian type. The precise statement and the proof are given in Section \ref{sec:asymp}. The result for the $p$-Laplace equation was due to \cite{kaminv88} and improved in \cite{LeePV06}.

In Section \ref{sec:bounded-dom} we  consider the problem posed in a bounded domain with zero boundary values. Again using an explicit solution that we  call a friendly giant type solution, we show that  the viscosity solution asymptotically behaves like such a solution, see the detailed statement in Theorem \ref{thm:bounded-dom}.

In  Section \ref{apriori}, we establish Harnack type estimates. We use the Barenblatt type solution for the expansion of positivity in the proof of Theorem \ref{thm:harnack}.  Combined with the oscillation estimate Corollary \ref{cor:key-osc}, this yields the  Harnack type estimate of Theorem \ref{thm:harnack}. The final  Appendix gathers a number of technical arguments.

\medskip

\noindent {\sc Notations and comments}. It is sometimes convenient to re-parametrize the general equation (\ref{eq:equation}) in the form
\begin{equation}\label{eq:qp}
\partial_t u=\abs{Du}^{\gamma}(\Delta u+(p-2)\Delta_{\infty}^{N} u)
=\abs{Du}^{q-2}(\Delta u+(p-2)\Delta_{\infty}^{N} u)\,,
\end{equation}
where $\gamma=\kappa +p-2>-1$ and $q=2+\gamma=\kappa +p>1$.  The value $\kappa=0$ ($\gamma=p-2$, $q=p$) indicates the standard  parabolic $p$-Laplace equation, while  $\gamma=0$  corresponds to the parabolic flow driven by the normalized $p$-Laplacian. We can also see the equation as a  particular  case of the general form \ $\partial_t u= F(Du,D^2 u)$ for a certain $F$, but we will not go far into this extra generality.

Finally, we point out that the problem becomes trivial in one space dimension. Indeed, since
then $\Delta_p u=(\abs{u_x}^{p-2}u_x)_x$, equation \eqref{eq:equation} reduces to a standard parabolic $q$-Laplace equation with the exponent $q=\kappa+p$:
$$
\partial_t u=\abs{u_x}^{\kappa}\Delta_p u= (p-1)/(q-1) \Delta_{q}u.
$$
 There is no need for a fictitious dimension in this case.

\section{Viscosity solutions}
\label{sec:visc}

Let $\Om\subset \Rn, n\ge 1$, be an open bounded set, and $\Theta \subset \R^{n+1}$ an open set. We define the cylinder $\Om_T=\Om\times(0,T)$ and its  parabolic boundary
\begin{equation}
\label{eq:parabolic-bounary}
\begin{split}
\partial_p \Om_T =( \Om\times \{0\})\cup (\partial \Om \times [0,T]).
\end{split}
\end{equation}
Moreover, we will use cylinders of the form $ Q_{r,\theta}(x_0,t_0)=B_{r}(x_0)\times (t_0-\theta, t_0),$ so that in particular $Q_{r,r^q}(x_0,t_0)=B_{r}(x_0)\times (t_0-r^q, t_0)$. When no confusion arises, we may drop the reference point and write $Q_{r,r^q}$. For a cylinder in  one space dimension, we denote $Q^1_{\delta,\delta}(x_0,t_0)=(\abs{x_0}-\delta,\abs{x_0}+\delta)\times (-\delta+t_0,t_0)$.

The definition of suitable viscosity solutions to (\ref{eq:qp}) requires some care because the operator may be singular. Nonetheless, a definition that fits our needs can be found in \cite{ohnumas97}. First set
\begin{equation}
\label{eq:F-def}
F(Du,D^2 u):=\abs{Du}^{q-2}(\Delta u+(p-2)\Delta_{\infty}^{N} u)
\end{equation}
whenever $Du\neq 0$. We define $\mathcal F$ to be a set of functions $f\in C^{2}([0,\infty))$ such that
\[
\begin{split}
 f(0)=f'(0)=f''(0)=0,\ f''(r)>0 \text{ for all }r>0,
\end{split}
\]
and moreover we require for  $g(x):=f(\abs x)$ that
\[
\begin{split}
\lim_{x\to 0,x\neq 0}F(Dg(x),D^2g(x))=0.
\end{split}
\]
Further, let
\[
\begin{split}
\Sigma=\{ \sigma \in C^{1}(\R) \,:\, \sigma \text{ is even},\, \sigma(0)=\sigma'(0)=0,\text{ and }\sigma(r)>0 \text{ for all } r\neq 0 \}.
\end{split}
\]

\begin{definition}
\label{def:admissible}
A function $\vp \in C^{2}(\Theta)$ is admissible if for any $(x_0,t_0)\in \Theta$ with $D\vp(x_0,t_0)=0$, there are $\delta>0$, $f\in \mathcal F$ and $\sigma \in \Sigma$ such that
\[
\abs{\vp(x,t)-\vp(x_0,t_0)-\vp_t(x_0,t_0)(t-t_0)}\le f(\abs{x-x_0})+\sigma(t-t_0)
\]
for all $(x,t)\in B_{\delta}(x_0)\times (t_0-\delta,t_0+\delta)$.
\end{definition}
If $D\vp\neq 0$, a $C^2$-function is automatically admissible.

\begin{definition}
\label{eq:from-below}
We say that $\vp$ touches $u$ at $(x_0,t_0)\in \Theta$ from below
if
\begin{enumerate}
\item[(1)] $u(x_0,t_0)=\vp(x_0,t_0)$, and
\item[(2)] $u(x,t)>\vp(x,t)$ for all $(x,t)\in \Theta$ such that $(x,t)\neq (x_0,t_0)$.
\end{enumerate}
\end{definition}

The definition for touching from above is analogous.

\begin{definition}
\label{def:visc}
A function $u:\Theta\to \R\cup \{\infty\}$ is a viscosity supersolution  to (\ref{eq:qp}) if
\begin{enumerate}
\item[(i)] $u$ is lower semicontinuous,
\item [(ii)] $u$ is finite in a dense subset of $\Theta$,
\item [(iii)]for all admissible $\vp\in C^{2}(\Theta)$ touching $u$ at $(x_0,t_0)\in \Theta$ \ from below
\[
\begin{split}
\begin{cases}
\vp_t(x_0,t_0)-F(D\vp(x_0,t_0),D^2\vp(x_0,t_0))\ge 0 & \text{if }D\vp(x_0,t_0)\neq 0,\\
\vp_t(x_0,t_0)\ge 0 &  \text{if }D\vp(x_0,t_0)= 0.
\end{cases}
\end{split}
\]
\end{enumerate}
\end{definition}

The definition of a subsolution $u: \Theta\to \R\cup \{-\infty\}$ is analogous except that we require upper semicontinuity, touching from above, and we reverse the inequalities above: in other words if $-u$ is a viscosity supersolution.  If a continuous function is both a viscosity super- and subsolution, it is a {\sl viscosity solution.}

It is shown in \cite{juutinenlm01} that if $q= p>1$, then the above notion coincides with the notion of $p$-super/subparabolic functions, having a direct connection to the distributional weak super/subsolutions as well. Moreover, if $q\ge 2$, then viscosity solutions can be defined in a standard way by using semicontinuous extensions, see Proposition 2.2.8 in \cite{giga06}.

The following comparison principle is proved in \cite[Theorem 3.1]{ohnumas97}.

\begin{theorem}[Comparison]
\label{thm:comparison}
Let $\Om$ be a bounded domain. Suppose that
 $u$ is viscosity supersolution and $v$ is a viscosity subsolution in $\Om_T$. If
 \[
\begin{split}
   \infty \ne    \limsup_{\Om_T \ni (y,s)\rightarrow (x,t)} v(y,s)\leq
    \liminf_{\Om_T \ni (y,s)\rightarrow (x,t)} u(y,s) \ne -\infty
\end{split}
\]
for all $(x,t) \in \partial_p\Om_T$,
then $v\leq u$ in $\Om_T$.
\end{theorem}
Actually, the result is proved  for the viscosity solutions of  the more general form $\partial_t u= F(Du,D^2 u)$ where $F$ satisfies certain regularity and degenerate ellipticity conditions, see  \cite[Section 2]{ohnumas97}.  They apply in our case for $p> 1$ and $q>1$.

This then implies the existence and uniqueness for the Cauchy problem for our problem (\ref{eq:qp}) by the Perron method \cite[Theorem 4.9]{ohnumas97}. Below $BUC$ refers to the space of bounded and uniformly continuous functions.

\begin{theorem}[Cauchy problem]
\label{thm:Cauchy-exist-uniq}
Let $u_0 \in BUC(\Rn)$. Then there exists a unique viscosity solution $u$ in the class $BUC(\Rn\times [0,T))$ to the problem
\begin{align*}
\begin{cases}
u_t=F(Du,D^2u), & \text{in }\Rn\times (0,T)\\
 u(x,0)=u_0(x),& x\in \Rn.
\end{cases}
\end{align*}
\end{theorem}

For a bounded open $\Om\subset \Rn$ with suitable regularity conditions a slight modification
of the Perron method also gives the existence of a unique viscosity solution to the Dirichlet problem with continuous boundary values $g$ (see \cite[Theorem 2.4.9]{giga06}). Another approach is to approximate the problem with a smoother one and to prove the existence then by passing to the limit. To be more precise, combining Theorem 5.2 and 5.3, Lemma 5.4, Theorem 5.5 in \cite{imbertjs}, the following theorem holds for $Q_T=B_1\times (0,T)$.

\begin{theorem}[Dirichlet problem]
Let $g\in  C(\partial_p Q_T)$. Then there exists a unique viscosity solution $u\in C(\overline{Q_T})$ to equation (\ref{eq:qp}) posed in $Q_T$ such that $u=g$ on $\partial_p Q_T$.
\end{theorem}

These solutions satisfy the $C^{1,\a}$ interior regularity for equation (\ref{eq:qp}) established  by Imbert, Jin and Silvestre in \cite{imbertjs}. The H\"older norm depends on the $L^\infty$ bound of the solutions and the domains.

\section{Radial solutions. Reduction to fictitious dimension}
\label{sec:radial}

In this section, we derive explicit radial solutions to  (\ref{eq:qp}) by a functional transformation. For more clarity, we introduce the transformation using formal computations. At the end of the section a theory of radial solutions is done using this transformation and we verify that the obtained formal solutions are indeed viscosity solutions.

\subsection{Introducing the problem in the fictitious dimension}
 Let  $u$ be a smooth solution  to (\ref{eq:qp}) which is radial with respect to the space variable, and has a non-vanishing gradient for $x\ne 0$. With a slight abuse of notation,  as  usual in the literature, we also use the same notation $u$ to denote the solution in radial coordinates, i.e. $u(x,t)$ becomes $u(r,t)$. Then, denoting by $u_r$ the radial derivative, we have $\abs{Du(x,t)}=\abs{u_{r}(r,t)}$, and
$$
\Delta u=u_{rr}+\frac{n-1}{r}u_{r}.
 $$
   Thus
\begin{equation}
\label{eq:radial-eq}
\begin{split}
u_t&=\abs{Du}^{q-2}(\Delta u+(p-2)\Delta_{\infty}^{N} u)\\
&=\abs{u_{r}}^{q-2}(u_{rr}+\frac{n-1}{r}u_{r}+(p-2)u_{rr})\\
&=\abs{u_{r}}^{q-2}(p-1)(u_{rr}+\frac{n-1}{(p-1)r}u_{r})\\
&=\frac{p-1}{q-1}\abs{u_{r}}^{q-2}\Big((q-1)u_{rr}+\frac{(n-1)(q-1)}{(p-1)r}u_{r}\Big).
\end{split}
\end{equation}
On the other hand, we can write the usual $q$-Laplacian for the smooth radial function  with a non-vanishing gradient in space dimension $d$ as
\begin{equation}
\label{eq:radial-p-Laplace}
\begin{split}
\Delta_{q}^d  u&:=\operatorname{div}(\abs{Du}^{q-2}Du)\\
&=\abs{Du}^{q-2}(\Delta u+(q-2)\Delta^N_{\infty}u)\\
&=\abs{u_{r}}^{q-2}\Big((q-1) u_{rr}+\frac{d-1}{r} u_{r}\Big),
\end{split}
\end{equation}
where we recall that $\Delta_{\infty}^{N} u=\sum_{i,j} \frac{\partial_{i}u}{|Du|}\frac{\partial_{j}u}{|Du|}\partial_{ij}u$ denotes the normalized infinity Laplacian.
Note that as long as we restrict to functions $u(r,t)$, the equations make sense even if $n$ and $d$ are not integers.  Comparing both formulas,  and starting from equation \eqref{eq:radial-eq}, we define the equivalent {\sl fictitious dimension $d(n,p,q)$}  as
\begin{equation}
\label{art.dim}
d=\frac{(n-1)(q-1)}{p-1}+1=\frac{(q-1)n+p-q}{p-1}\,.
\end{equation}
We may then write the radial equation in the form
\begin{equation}
\label{eq:1D-div-form}
\begin{split}
u_t&=\frac{p-1}{q-1}\abs{u_{r}}^{q-2}\Big((q-1) u_{rr}+\frac{d-1}{r} u_{r}\Big)\,,
\end{split}
\end{equation}
which is formally the parabolic $q$-Laplace equation in $d$ dimensions with the constant, an equivalence that will lead to interesting conclusions. Some remarks first:

(i) This type of transformation was used in \cite{portilheirov12, portilheirov13} when $p=\infty$, but then $d=1$ which makes things much easier. This is never the case here and we get $d>1$  whenever $p,q>1$ and $n>1$.

(ii) We have
$$
d-n= \frac{(n-1)(q-p)}{p-1}.
$$
Therefore, when  $q>p$  the fictitious dimension is larger than $n$, if $q<p$ then $d$ is less than $n$, and if $p=q$ there is naturally no change.
Conversely, given $d,n$ and $p$ we can get the needed exponent change as
$q-p=(p-1)(d-n)/(n-1). $

(iii)  The fictitious dimension may or may not be an integer. If $d$ is an integer, then (\ref{eq:qp})  and (\ref{eq:1D-div-form}) formally coincide  in the case of the radial solutions with the $q$-parabolic equation in $d$ space dimensions $u_t=\Delta^{d}_q u$ as shown by the above computations. The equivalence of (\ref{eq:qp}) and the 1-dimensional equation (\ref{eq:1D-div-form}) in the case of the radial solutions is discussed in Theorem \ref{thm:equiv}.

(iii) If $d$ is not an integer, the transformation still implies that $u(r,t)$ satisfies a 1-dimensional parabolic equation. This will be used below to construct examples that will be quite useful in the sequel.

\subsection{Examples}
\label{ss.examples}

The transformation we have introduced is useful because it produces interesting examples that we need later in the paper to settle the long-time behaviour of general solutions, and to prove Harnack's inequality.

\noindent {\bf Example 1}.
 We consider the simplest case $q=2$ where the right-hand operator has linear homogeneity, and look for source-type solutions, i.\,e., solutions that start from initial data consisting of a singularity. The fictitious dimension is now
$$
 d=\frac{p+n-2}{p-1}>1,
$$
and we see that (\ref{eq:radial-eq}) just becomes the heat equation in (formal) dimension $d$, but with the  factor $p-1$ in the equation, which is absorbed into the time variable. From the explicit fundamental solution of the heat equation, we derive  the formal solutions
\begin{equation} \label{eq:heat-type}
U(x,t)= C\,t^{-\frac{d}{2}}\exp\Big(-\frac{\abs{x}^2}{4(p-1)t}\Big)=
C\,t^{-\frac{n+p-2}{2(p-1)}}\exp\Big(-\frac{\abs{x}^2}{4(p-1)t}\Big)\,.
\end{equation}
These are the solutions obtained in \cite{banerjeeg13} to $u_t=\Delta u+(p-2)\Delta_{\infty}^N u$. In case $d$ is not an integer this is only formal, but it  can be shown  that they are viscosity solutions, see Proposition \ref{prop:barenblatt-is-visc}.

\medskip

\noindent {\bf Example 2}. For the other values $q\neq 2$, we get a different source-type solution that comes from the Barenblatt-type solution \cite{barenblatt52} of the standard parabolic $q$-Laplace equation in $d$ space dimensions. We will call the obtained solution a {\sl modified Barenblatt solution}. It is convenient to consider first the case  $q>2$ corresponding to slow diffusion. The standard Barenblatt solution reads then
\begin{equation}
\label{eq:barenblatt}
\begin{split}
  {\mathcal B}^{q,d}(x,t;C)=t^{-d/\lambda}\Bigg(C-\frac{q-2}{q}\lambda^{\frac{1}{1-q}}\Big(\frac{\abs x}{t^{1/\lambda}}\Big)^{q/(q-1)}\Bigg)_+^{\frac{q-1}{q-2}}\,,
\end{split}
\end{equation}
where $\lambda=d(q-2)+ q$,  the constant $C>0$ can be chosen freely, and $(\cdot)_+$ means $\max\{\cdot,0\}$. Notice the property of compact support, and this property only depends on the condition $q>2$.

Using the value \eqref{art.dim} and doing the time rescaling we arrive at a formal expression for the source type solution of (\ref{eq:qp}) of the form
\begin{equation}
\label{eq:barenblatt-q-d}
\begin{split}
 &\mathcal B^{q,d,p}(r,t;C):=\mathcal B^{q,d}(re_1, \frac{p-1}{q-1}t;C) \text{ and }\\
   &\mathcal B^{q,d,p}(x,t;C):=\mathcal B^{q,d,p}(\abs{x},t;C)
\end{split}
\end{equation}
where $e_1=(1,0,\ldots,0)$. Sometimes  we drop $C$, i.\,e., \ denote $  \mathcal B^{q,d,p}(r,t)$ and $ \mathcal B^{q,d,p}(x,t)$ for simplicity.

 The limit $p\to\infty$ can be taken in these examples, and it leads to the results of \cite{portilheirov12, portilheirov13} with $d=1$.

\medskip

\noindent {\bf Example 3}.  On the other hand, for  $1<q<2$ we have a fast diffusion Barenblatt solution if $ \lambda $ remains positive, i.e., for $d(q-2)+q>0$, see \cite{vazquez06} , page 192 or \cite{dibenedettogv11}. Because of the sign change the formula is now
\begin{align*}
 \mathcal B^{q,d}(x,t)=t^{-d/\lambda}\Bigg(C+\frac{2-q}{q}\lambda^{\frac{1}{1-q}}\Big(\frac{\abs x}{t^{1/\lambda}}\Big)^{q/(q-1)}\Bigg)^{-\frac{q-1}{2-q}}\,,
\end{align*}
which defines the time rescaled version
\begin{equation}
\label{eq:barenblatt.fast}
 \mathcal B^{q,d,p}(x,t;C)
\end{equation}
similarly as in (\ref{eq:barenblatt-q-d}).
The outcome in \eqref{eq:barenblatt-q-d} is the same, but in this example the solutions do not have compact support, they have instead a tail with power-rate decay as $\abs x\to\infty$.

Let us examine the admissible range for $p$ and $q$ when $q<2$. We have to impose the condition
$$
\lambda =d(q-2)+q=\Big(\frac{(n-1)(q-1)}{p-1}+1\Big)(q-2)+q>0\,,
$$
i.e.
$$q>\frac{2d}{d+1}$$
which is equivalent to $(n-1)(2-q)<2(p-1)$, i.\,e.,
\begin{equation}
\label{eq:pre-range}
\begin{split}
 2n < q(n-1)+2p.
\end{split}
\end{equation}
It follows that the range condition can be written as
\begin{equation}
\label{eq:range}
\begin{split}
q>\begin{cases}
1 & \text{ if }p\ge (1 + n)/2 \\
2 (n - p)/(n-1) &\text{ if } 1<p< (1 + n)/2 .
\end{cases}
\end{split}
\end{equation}
Observe that $1<(2 (n - p))/(n-1)<2$ whenever $1<p< (1 + n)/2$.

\medskip

\noindent {\bf Conservation laws.} (i) \kom{naturally: $q>2 d/(d+1)$} The change of dimension when $q\ne p$ does not seem important in the given formulas, but it has consequences for the physical interpretation as we will see below. Let us point out a very important fact: the standard Barenblatt solution in dimension $d$ obeys the mass conservation with respect to the Lebesgue measure   $dx=\ud S\ud r$, so that in radial coordinates we have
\begin{equation}
\int_0^\infty  \mathcal B^{q,d,p}(r,t)r^{d-1}dr= C,
\end{equation}
where $C>0$ is a constant independent of time.  We call this integral the $d$-mass. Note that this fact is a consequence of the self-similar form of the Barenblatt solutions with the self-similarity exponents given below.

 The equation \eqref{eq:qp} is invariant under a scaling transformation  of the form
\begin{equation}\label{scaling}
{\mathcal T}u(x,t)=A\,u(B x, C t)
\end{equation}
for real parameters $A,B,C>0$. This formula transforms solutions into solutions if $C=A^{q-2}B^{q}$. This leaves two free parameters
which can be conveniently used in the theory.

 The conservation of the $d$-mass above is a consequence of the fact that the Barenblatt solutions  obey the above type scaling invariance  with the extra $d$-mass condition
$A=B^d$, so that $C=B^{(q-2)d+q}$ i.e.\ $u(x,t)=B^d u(Bx,B^{\lambda} t)$.

(ii) When we try to write this conservation law for the modified Barenblatt solution with respect to the $n$-dimensional measure we get
\begin{equation}
\int_0^\infty r^{\sigma}\,  {\mathcal B}^{q,d,p}(r,t)\,r^{n-1}dr= C,
\end{equation}
which means conservation of the moment taken with respect to the standard Lebesgue measure  with the weight $w(r)=r^{\sigma}$ where
$$
\sigma = d-n= \frac{(n-1)(q-p)}{p-1}.
$$

\noindent {\bf Initial singularity.}  In the same vein, the calculation for the $n$-mass gives
$$
\int_0^\infty \, {  \mathcal B^{q,d,p}(r,t)}\,r^{n-1}dr= C\,t^{-\mu}, \quad
 \mu=\frac{d-n}{\lambda}=\frac{(n-1)(p-q)}{(p-1)(d(q-2)+q)}\,.
$$
In view of the fact that the Barenblatt solution for $q>2$ has a support that shrinks to the origin as $t\to 0$, we conclude that the initial data is a Dirac delta only if $p=q$. If $q>p$ the initial mass tends to infinity (infinite mass singularity), if $q<p$ it tends to zero (a mild singularity).

\medskip

\noindent {\bf Justification of the formulas}.  The derivation of the Barenblatt solutions is well known for integer space dimensions, and they are  weak solutions of the equation. We need to justify that such claims hold when $d$ is not an integer.

Here are the whole details   of the formal pointwise computation of the solution.
 We look for a nonnegative radial self-similar solution. We propose the form
\[
\begin{split}
v(r,t)=(\beta t)^{-\frac{\a}{\beta}}w((\beta t)^{-\frac1{\beta}}r)\,,
\end{split}
\]
where $\alpha$ and $\beta$ are to be suitably fixed later,  $r, t>0$, and we assume that all the derivatives below  exist. We get
\[
\begin{split}
v_t(r,t)&=-\frac{\a}{\beta} (\beta t)^{-\frac{\a}{\beta}-1}\beta w((\beta t)^{-\frac1{\beta}}r)-  \frac1{\beta}(\beta t)^{-\frac{\a}{\beta}}(\beta t)^{-\frac1{\beta}-1}r w'((\beta t)^{-\frac1{\beta}}r)\beta\\
&=-(\beta t)^{-\frac{\a+\beta}{\beta}}(\a w(R)+Rw'(R)),
\end{split}
\]
where $R:=(\beta t)^{-\frac1{\beta}}r$.  Moreover,
\[
v_{r}(r,t)=(\beta t)^{-\frac{\a+1}{\beta}}w'(R), \quad
v_{rr}(r,t)=(\beta t)^{-\frac{\a+2}{\beta}}w''(R).
\]
Inserting this into $v_t=\abs{v_{r}}^{q-2}((q-1)v_{rr}+\frac{d-1}{r}v_{r})$ we obtain
\[
\begin{split}
-(\beta t)^{-\frac{\a+\beta}{\beta}}&(\a w(R)+Rw'(R))\\
&=\abs{(\beta t)^{-\frac{\a+1}{\beta}}w'(R)}^{q-2}\Big((q-1)(\beta t)^{-\frac{\a+2}{\beta}}w''(R)+\frac{d-1}{r}(\beta t)^{-\frac{\a+1}{\beta}}w'(R)\Big)\\
&=(\beta t)^{-\frac{\a+2}{\beta}}(\beta t)^{-\frac{(\a+1)(q-2)}{\beta}}\abs{w'(R)}^{q-2}\Big((q-1)w''(R)+\frac{d-1}{R}w'(R)\Big).
\end{split}
\]
We may eliminate the time dependence by choosing $\a+\beta=(\a+2)+(\a+1)(q-2)$ i.\,e., \  $\beta=\a(q-2)+q>0$, and  we get
\[
\begin{split}
-(\a w(R)+Rw'(R))=\abs{w'(R)}^{q-2}\Big((q-1)w''(R)+\frac{d-1}{R}w'(R)\Big).
\end{split}
\]
Moreover,
\begin{equation}
\label{eq:tb-integrated}
\begin{split}
-R^{1-\a}(R^{\a}w(R))'= \Big(\abs{w'(R)}^{q-2}w'(R)R^{d-1}\Big)' R^{1-d},
\end{split}
\end{equation}
where $(\cdot)'$ denotes the derivative with respect to $R$ and we used
\begin{equation}
\label{eq:1-divergence-form}
\begin{split}
 &\Big(\abs{w'(R)}^{q-2}w'(R)R^{d-1}\Big)' R^{1-d}\\
 &=\Big(\sgn(w'(R))(q-2)\abs{w'(R)}^{q-3}w''(R)w'(R)R^{d-1}+\abs{w'(R)}^{q-2}w''(R)R^{d-1}\\
 &\hspace{10 em}+\abs{w'(R)}^{q-2}w'(R)(d-1)R^{d-2}\Big) R^{1-d}\\
&=\abs{w'(R)}^{q-2}\Big((q-1)w''(R)+w'(R)\frac{d-1}{R}\Big)
\end{split}
\end{equation}
assuming $w'(R)\neq 0$; the vanishing gradient will need a special treatment later. For the later use, we also remark that if $w$ is extended through an even reflection $w(R):=w(-R)$, $R<0$, then the above equation takes the form
\begin{align*}
\Big(\abs{w'(R)}^{q-2}w'(R)\abs{R}^{d-1}\Big)'\abs{R}^{1-d}=\abs{w'(R)}^{q-2}\Big((q-1)w''(R)+w'(R)\frac{d-1}{R}\Big).
\end{align*}

The choice $\a=d$ in (\ref{eq:tb-integrated}) together with an integration gives
\begin{equation}
\label{eq:integrated}
\begin{split}
0=Rw(R)+ \abs{w'(R)}^{q-2}w'(R).
 \end{split}
\end{equation}
A solution to this reads as
\[
\begin{split}
w(R)=\Big(K-\frac{q-2}{q}R^{\frac{q}{q-1}}\Big)_{+}^{\frac{q-1}{q-2}}
\end{split}
\]
and thus
\[
\begin{split}
v(r,t)&=(\beta t)^{-\frac{\a}{\beta}}w((\beta t)^{-\frac1{\beta}}r)=t^{-d/\beta}\Big(C-\frac{q-2}{q}\beta^{\frac1{1-q}}\Big(\frac{r}{t^{1/\beta}}\Big)^{\frac{q}{q-1}}\Big)_{+}^{\frac{q-1}{q-2}}\\
&=t^{-d/\lambda}\Big(C-\frac{q-2}{q}\lambda^{\frac1{1-q}}\Big(\frac{r}{t^{1/\beta}}\Big)^{\frac{q}{q-1}}\Big)_{+}^{\frac{q-1}{q-2}},
\end{split}
\]
where $C=\lambda^{-d(q-2)/(\lambda(q-1))}K$ and we recalled the earlier notation $\lambda=\beta$ from (\ref{eq:barenblatt}). Finally, letting $u(x,t):=v(x,\frac{p-1}{q-1}t)$ solves, at least formally at this point, the equation (\ref{eq:1D-div-form}).

 \medskip

The presence of a free boundary where the regularity is limited for $q>2$ implies that a proof is needed to show that it is indeed a viscosity solution of the equation.  In Proposition \ref{prop:barenblatt-is-visc} below we show that this is the case  for all values of $d>1$.

\medskip

\noindent {\bf Other examples.} Similar procedures can be applied to other families of explicit solutions of the standard $p$-Laplace equation to produce new solutions of equation \eqref{eq:qp}.  For instance, a simple solution that is sometimes used as a barrier is
$$
u(x,t)=c(at+x_1-b)_+^{(p-1)/(p-2)}
$$
with  $p>2$, an arbitrary $a>0$ and a convenient $c=c(a,p)$, and any $b\in\re$. This is a one-dimensional traveling wave directed along the $x_1$-axis (for any other direction we get a solution by rotation of this one).
In this case, the corresponding solution for equation \eqref{eq:qp} is given by the same formula with  $p$ replaced by $q$, and we do not need any change of dimension. We will not exploit the last example further in this paper.


\section{Equation in 1-D. Basic results}\label{sec.1dim}

The  discussion of the last section motivates the study of the 1-dimensional parabolic equation  \eqref{eq:1D-div-form} i.e.
\begin{equation}\label{eq.qd}
u_t=\frac{p-1}{q-1}\abs{u_{r}}^{q-2}\Big((q-1) u_{rr}+\frac{d-1}{r} u_{r}\Big).
\end{equation}

\subsection{Recapitulation of the theory for integer $d$}

Whenever $d$ is an integer, the equation \eqref{eq.qd} is just the standard $q$-Laplace equation in $\re^d$ for radial functions. Radial viscosity solutions of the original equation transform into radial viscosity solutions of the $d$-dimensional $q$-Laplace equation,  cf.\ the estimate \eqref{eq:CP-from-below}.

The theory of the standard  $d$-dimensional $q$-Laplace equation is well-established and for example the following holds:

(i) Given data $u_0\in L^{ s }(\Om)$,  where $s\in [1, \infty)$ is any exponent   there is a unique weak solution of both problems ($\Om$ bounded or $\R^d$), the set of solutions forms a contraction semigroup in $L^{s}(\Om)$, see  \cite{Barbu, BenCr91, Br73, dibenedetto93, Rv70}.

(ii) For bounded initial data, the solutions  $u$ are locally  $C^{\alpha}$ in space and time for some $\alpha\in (0,1)$.  Moreover, $D u\in C_{x,t}^{\alpha,\alpha/2}$,  \cite{dibenedetto93}. Continuous data produce continuous solutions up to the boundary in regular domains \cite{dibenedetto93, kilpelainenl96, bjornbgp15}.
For $q\ge 2$ all $L^{s }(\Om)$-solutions, $s\ge1$, are bounded; this is also true for $1<q<2$ if  $q$ is not too small, $q(d+1)> 2d$. However, we emphasize that viscosity solutions are always continuous by definition.

(iii) For $q>2$ we have {\sl finite speed of propagation}. We formulate it in the simplest form: if $\Om=\re^d$ and the initial data are nonnegative, bounded and supported in a finite ball $B_{R_0}(0)$, then for all $t>0$ the support of the solution $u(\cdot,t)$ is contained in a finite ball of finite radius $R(t)$ and $R(t)\to \infty$ as $t\to \infty$. The support of the solution increases with time.

(iv) On the contrary, for $1<q<2$ we have {\sl infinite speed of propagation}: every nonnegative, continuous and bounded weak solution of the stated problems defined in $\Om\times (0,T)$ will be either positive everywhere or identically zero for each $0<t<T$. Actually, there is a time $T_e>0$ such that $u(x,t)>0$ for every $x\in\Om$ and $0<t<T_e$ and $u(x,t)=0$ $x\in\Om$ and $T_e< t\le T$.
When $T_e<  T$ then $T_e$  is called the extinction time. It depends on the problem and on the initial data. Note that for the Cauchy-Dirichlet problem in a bounded domain with zero lateral boundary data we may take $T= \infty$ and $T_e$ is always finite.
On the other hand, for the Cauchy Problem in $\re^d$ we may take $T=\infty$ (global solutions) and finite time extinction depends on $q$ and the class of data. Thus,  if $u_0\in L^1(\re^d)$, then $T_e(u_0)$ is finite whenever  $1<q<2d/(d+1)$, while for  $2>q\ge 2d/(d+1)$ the conservation of mass holds. See more for example \cite{vazquez06} (Chapter 11) or \cite{dibenedetto93} (Chapter 7).

\subsection{The case of non-integer $d$}

Here we  will establish the basic facts that will allow us to work with the transformations of the viscosity solutions of equation \eqref{eq:qp}. We will choose to present the details for the Cauchy-Dirichlet problem with continuous initial data.
 Below use the notation \ $\ud z:=\abs{r}^{d-1}\ud r \ud t$,  the natural parabolic measure for this problem.

\begin{definition}
\label{def:1D-weak}
Let $0<T\le \infty$ and $0<R\le \infty$.
A function $u \in C((-R,R)\times (0,T))$ such that  $u_r \in C((-R,R)\times (0,T)), u_r(0,t)=0$ is a continuous weak solution to (\ref{eq:1D-div-form}) if
\[
\begin{split}
\int_{(-R,R)\times(0,T)} u \phi_t \ud z  = \frac{p-1}{q-1} \int_{(-R,R)\times(0,T)} \abs{u_{r}}^{q-2}u_{r} \phi_{r} \ud z
\end{split}
\]
for all $\phi\in C^\infty_0((-R,R)\times(0,T))$.  To get the definition of a continuous  weak subsolution, we replace equality by $\ge$ and test with $\phi\ge 0$. The definition for a continuous weak supersolution is analogous except the inequality is reversed.
\end{definition}

Above we have taken the rather strong regularity assumption for convenience, since the corresponding viscosity solutions are even in $C^{1,\alpha}$ and this will be the context  where we use the definition. Let us also remind that in the equivalence theorems below we consider radial, and in the asymptotic results in the whole $\Rn$, we assume boundedness.

Next we observe that the  radial viscosity solutions of Section \ref{sec:visc}  and the 1-dimensional continuous weak solutions we have introduced are the same. Observe that since we assume in the next theorem that the function is radial and necessarily $u_r(0,t)=0$, then we have even function with respect to $r$.
 \begin{theorem} \label{thm:equiv}
Let $u\in C(Q_T)$, $Q_T=B_R\times(0,T), B_R\subset \R^n, 0<R\le \infty, $ be a continuous radial function, and $q>1$.
Then $u$ is a viscosity solution to (\ref{eq:qp}) in $n$-dimensions if and only if  $v(r,t):=u(re_1,t),\ r\in (-R,R),$ is 1-dimensional  weak solution to (\ref{eq:1D-div-form})
according to Definition \ref{def:1D-weak}.
\end{theorem}

We have decided to postpone the proof of this result to
 the appendix, see Propositions \ref{prop:weak-is-visc} and \ref{prop:visc-is-weak}, to give precedence to the asymptotic results in Section \ref{sec:asymp}.
\begin{remark}
\label{rem:extensions}
(i) In case $d$ is an integer, similarly it holds that the radial viscosity solutions in $n$-dimensions to (\ref{eq:qp}) i.e.
\begin{align*}
u_t=\abs{Du}^{q-2}(\Delta u+(p-2)\Delta_{\infty}^{N} u)
\end{align*}
are equivalent to the radial weak solutions to the parabolic $q$-Laplacian in $d$-dimensions
\begin{align*}
u_t&=\frac{p-1}{q-1}\Delta_{q}^d u.
\end{align*}
(ii) A similar equivalence result also holds in the time independent case: Let $u\in C(B_R)$,  $B_R\subset \R^n, 0<R< \infty, $ be a continuous radial function, $q,p>1$, and  $f$ radial as well as continuous up to the boundary.
Then $u$ is a viscosity solution to 
\begin{align}
\label{eq:elliptic}
\abs{Du}^{q-2}(\Delta u+(p-2)\Delta^N_{\infty} u)=f
\end{align}
 in $n$-dimensions if and only if  $v(r):=u(re_1),\ r\in (-R,R),$ is 1-dimensional  weak solution to 
 \[
\begin{split}
\int_{(-R,R)} f \phi \abs{r}^{d-1} \ud r  = -\frac{p-1}{q-1} \int_{(-R,R)} \abs{u_{r}}^{q-2}u_{r} \phi_{r} \abs{r}^{d-1}\ud r,
\end{split}
\]
for all $\phi\in C^\infty_0((-R,R))$.
If $d$ happens to be an integer, then in the radial case the viscosity solutions are equivalent with the weak solutions to the equation
\begin{align*}
\frac{p-1}{q-1} \Delta^{d}_{q} u=f
\end{align*}
where $ \Delta^{d}_{q}$ denotes the standard $q$-Laplacian in $d$-dimensions. 

Above we require for example  $u \in C((-R,R)),\ u_r \in C((-R,R)),\ u_r(0,t)=0$ in the weak definition. The $C^{1,\a}$-regularity for viscosity solutions of (\ref{eq:elliptic}) was proven in \cite{birindellid12} in the radial case; for the general case see \cite{attouchir}.
\end{remark}
\medskip

We now state a very remarkable property of this equation, in line with what was said for the Barenblatt solutions.

\noindent {\bf Conservation law}.  By Theorem \ref{thm:equiv}  radial viscosity solutions of the original equation  \eqref{eq:qp} coincide with the 1-dimensional solutions of the divergence form equation with weights and thus have a conservation of mass property in the fictitious dimension $d$, i.e.,
\[
\begin{split}
 \int_{0}^{\infty} u(r,t)r^{d-1}\ud r
\end{split}
\]
is constant in time for all $t> 0$ under a suitable global condition. It is well known that the conservation of mass fails even for the heat equation without a global condition, e.g.\ exponential growth bound.
 Here,  $d$ need not be an integer. We do not know of any conservation law of this type for general nonradial viscosity solutions.

The conservation law also holds in the singular range $q<2$ as long as the range condition (\ref{eq:range}) i.e.\ $q>2 d/(d+1)$ and a suitable global condition holds. Indeed, the proof in \cite{finodfv14} also holds  for our 1-dimensional equation. For the standard $q$-Laplacian $q>2$, see \cite{dibenedetto93}, Chapter 7.

\subsection{Back to the Barenblatt solutions}

From Theorem \ref{thm:equiv} it follows that the explicit solutions given in Subsection \ref{ss.examples} are viscosity solutions.

\begin{proposition}
\label{prop:barenblatt-is-visc}
The function $\B$ in (\ref{eq:barenblatt-q-d}) is a viscosity solution to (\ref{eq:qp}) in $\Rn\times (0,\infty)$ for $q>2$, and so is (\ref{eq:barenblatt.fast})  whenever $q<2$ and the range condition (\ref{eq:range}) holds. If $q=2, p>1$, then a viscosity solution is given by (\ref{eq:heat-type}).
\end{proposition}
\noindent {\sc Proof.}
By Theorem \ref{thm:equiv}, it suffices to show that the solutions in the statement are weak solutions according to Definition \ref{def:1D-weak}. Also observe that the definition of viscosity solutions still applies since it is of local nature, and thus the initial singularity is not a problem.

{\bf Case $q=2$:}  First, since (\ref{eq:heat-type}) is a smooth classical solution to (\ref{eq:1D-div-form}), then this is immediate in the case  $q=2$.

{\bf Case $q\neq2$:} Set
 \[
\begin{split}
A(r,t):=\Bigg(C-\frac{q-2}{q}\lambda^{\frac{1}{1-q}}\Big(\frac{r}{t^{1/\lambda}}\Big)^{q/(q-1)}\Bigg).
\end{split}
\]
If $A(r,\frac{p-1}{q-1}t)\neq 0,\ t>0,\ r\neq 0$, then $\B$ is a classical solutions to (\ref{eq:1D-div-form}) by construction and thus weak, and if $r=0$ then $(\B)_r(r,t)=0$ and one can directly verify the weak definition at the vicinity of $r=0$.

 Moreover, if $A(r,\frac{p-1}{q-1}t)=0,\ t>0$, then $(\B)_t$ and $(\B)_r$ are continuous and one can  again directly verify the definition of the weak solution.

  By Theorem \ref{thm:equiv}, it suffices to show that  $(\ref{eq:barenblatt-q-d})$ or (\ref{eq:heat-type}) are weak solutions according to Definition \ref{def:1D-weak}.
 \qed


For completeness, we also give  a more direct proof using the definition of viscosity solutions.

\medskip

%
\noindent {\sc A second proof of Proposition \ref{prop:barenblatt-is-visc}.} We consider $(x_0,t_0)\in \Rn\times (0,\infty)$.

\noindent {\bf Case $q>2$:} Let
\[
\begin{split}
A(x,t):=\Bigg(C-\frac{q-2}{q}\lambda^{\frac{1}{1-q}}\Big(\frac{\abs x}{t^{1/\lambda}}\Big)^{q/(q-1)}\Bigg).
\end{split}
\]
If $x_0\neq 0$ and $A(x_0,\frac{p-1}{q-1}t_0)\neq 0$, then $\B$ immediately satisfies the definition of a viscosity solution by the calculations in the examples of Subsection \ref{ss.examples}, since either $D\B(x_0,t_0)\neq 0$ or $\B$ is identically zero in the neighborhood of the point $(x_0,t_0)$. Indeed, let $\vp$ touch $\B$ from below at $(x_0,t_0)$ where $D\B(x_0,t_0)\neq 0$. Then observe that $D\vp(x_0,t_0)=D\B(x_0,t_0)$, and $D^2\vp(x_0,t_0)\le D^2 \B(x_0,t_0)$. Thus
\[
\begin{split}
F(D\vp,D^2\vp)&\le F(D\B,D^2\B)=\partial_t \B=\partial_t \vp.
\end{split}
\]

If $x_0=0$, then it holds that
\[
\begin{split}
\lim_{0\neq x\to 0} F(D\B(x,t_0),D^2\B(x,t_0))&=\partial_t\B(0,t_0) \\
&=-\frac{d}{\lambda} \Big(\frac{p-1}{q-1}t_0\Big)^{-\frac{d}{\lambda}-1}\frac{p-1}{q-1} C^{\frac{q-1}{q-2}}<0. 
\end{split}
\]
There is no admissible test function from below, cf.\ \cite{ohnumas97} Section 5. Hence the supersolution property is automatically satisfied. On the other hand,   since $\partial_t\B(0,t_0)<0$, then any admissible test function from above satisfies the definition of a viscosity subsolution.

If $A(x_0,\frac{p-1}{q-1}t_0)=0$ (i.e.\  we are at the free boundary), then since $\frac{q-1}{q-2}-1=\frac{1}{q-2}>0$, it follows that
\[
\begin{split}
D \B(x_0,t_0)=0=\partial_t \B(x_0,t_0).
\end{split}
\]
Thus for any test function $\vp$ touching $u$ at $(x_0,t_0)$ it holds that $\partial_t \vp(x_0,t_0)=0$ and thus the definition of a viscosity solution is satisfied.

\noindent {\bf Case $q<2$:} First we observe that in this case $A(x_0,\frac{p-1}{q-1}t_0)>0$, and $D\B(x_0,t_0)\neq 0$ if $x_0\neq 0$. In this case the Barenblatt solution  is also classical, and the above argument holds verbatim. Also the case $x_0=0$ can be treated exactly as above.

\smallskip

\noindent {\bf Case $q=2$:}  The solution (\ref{eq:heat-type})
\begin{align*}
U(x,t)=C\,t^{-\frac{d}{2}}\exp\Big(-\frac{\abs{x}^2}{4(p-1)t}\Big)
\end{align*}
is smooth and has nonzero gradient whenever $x_0\neq 0$.
In this case for any test function $\vp$ touching $U$ at $(x_0,t_0)$ from below
\begin{align*}
\partial_t \vp(x_0,t_0)=\partial_t U(x_0,t_0)=F( DU(x_0,t_0), D^2U(x_0,t_0))\ge F( D\vp(x_0,t_0), D^2\vp(x_0,t_0)).
\end{align*}
The case of touching from above is similar. If $x_0=0$, then again
\[
\begin{split}
\lim_{0\neq x\to 0} F(DU(x,t_0),D^2U(x,t_0))=\partial_tU(0,t_0) <0
\end{split}
\]
and there is no test function from below. When testing from above, we always have $
\partial_t\vp(0,t_0)=\partial_tU(0,t_0)<0$, and the definition of the viscosity solution is automatically satisfied.
\qed

\section{Asymptotic behaviour in the whole space}
\label{sec:asymp}

 We will now proceed with  the study of the long time behaviour. We will establish the asymptotic behaviour of the viscosity solutions of the general equation \eqref{eq:qp} in two typical situations: the initial-value problem in the whole space, and the Cauchy-Dirichlet problem with the zero lateral boundary data posed in a bounded domain. We point out that the techniques can be applied to other situations as well.

\subsection{Size estimates}

We consider first the initial-value problem to  equation \eqref{eq:qp} posed in $\ren$ and $F$ is given by (\ref{eq:F-def})
\begin{align}
\label{eq:Cauchy-problem}
\begin{cases}
u_t=F(Du,D^2u), & \text{in }\Rn\times (0,T)\\
 u(x,0)=u_0(x),& x\in  \Rn,\ 0\le u_0\in C_0(\Rn), u_0 \not\equiv 0
\end{cases}
\end{align}
with bounded and continuous solutions in the viscosity sense. $ C_0(\ren)$ denotes the space of continuous, compactly supported functions. By Theorem \ref{thm:Cauchy-exist-uniq} we know that this problem has a unique viscosity solution.

\medskip

 We want to obtain first a rough estimate of the size of the solutions, and also the free boundaries when $q>2$. A direct comparison with the explicit solutions constructed in Section \ref{sec:radial} produces a first bound on the solutions. These estimates are correct for all large times up to constant factors in view of the sharper results to follow.

\begin{theorem}
\label{thm:Cauchy-asymp.a}
Let $u$ be a viscosity solution to the Cauchy problem (\ref{eq:Cauchy-problem})  and the range condition (\ref{eq:range})
holds.

\noindent (i) There are constants $C,t_1>0$ such that
\begin{equation}
u(x,t)\le \mathcal B^{q,d,p}( x,t+t_1;C)    \qquad x\in \ren, \ t>0.
\end{equation}

\noindent (ii) There are constants $C_1,t_1>0$ such that
\begin{equation}
u(x,t)\ge \mathcal B^{q,d,p}( x,t+t_1;C_1)    \qquad x\in \ren, \ t>t_1>0.
\end{equation}

\noindent (iii) This means that for large times
\begin{equation}\label{alpha.f}
\|u(\cdot,t)\|_\infty \sim  \, t^{-\alpha}, \quad \alpha=\frac{d}{\lambda}=\frac{1}{q-2+(q/d)}\,,
\end{equation}
where $d=(n(q-1)+p-q)/(p-1)$ is the fictitious dimension, $\lambda=d(q-2)+q$,  $C_1$ depends on the initial data, and $\sim$ means up to a constant.

\noindent (iv) Moreover, when $q>2$ the support of the solution will be contained for large times in a ball of radius
\begin{equation}
R(t)= C_2\, t^{\alpha/d}\,, \quad \frac{\alpha}{d}=\frac{1}{d(q-2)+q}\,.
\end{equation}
while for $q<2$ the solution decays as $\abs x\to\infty$ like
\begin{equation}
\label{eq:tail}
u(x,t)\le C(t)|x|^{-q/(2-q)},
\end{equation}
where $C(t)$ is a given function of $t$.
\end{theorem}
\noindent {\sc Proof.}
\noindent {\bf Case $q>2$:}  Since $u_0\in  C_0(\Rn)$, $0\not\equiv u_0\ge 0$, we have Barenblatt type solutions as in (\ref{eq:barenblatt-q-d}) denoted by $\mathcal B_l(x,t):=\mathcal B^{q,d,p}(x-x_0,t+1; {C_l})$, $\mathcal B_u(x,t):=\mathcal B^{q,d,p}(x-x_0,t+1; {C_u})$ (lower, upper) such that
\[
\begin{split}
\mathcal B_l(x,0)\le u_0(x) \le \mathcal B_u(x,0).
\end{split}
\]
From the comparison principle, Theorem \ref{thm:comparison}, it follows for $(x,t)\in \R^n\times[0,\infty)$ that
\[
\begin{split}
\mathcal B_l(x,t)\le u(x,t) \le \mathcal B_u(x,t).
\end{split}
\]
This step fixes the size both from above and below with estimates that are proportional for large times.

\medskip

\noindent {\bf Case $q\le 2$:}   A modification has to be done to obtain the bounds from below since $\mathcal B_l(x,t)$ does not have compact support. Here is the argument in short. Take a ball $B_R(x_0)$ where $u_0$ is continuous and strictly positive, $u_0(x)\ge c>0$. By continuity of the solutions $u(x,t)\ge  c/2$ for  $x\in B_R(x_0)$ and $0<t<t_1$. We now consider the exterior space-time domain $E_{t_1}=(\ren\setminus B_R(x_0))\times (0,t_1)$.
We take a Barenblatt solution $\mathcal B_l(x,t):=\mathcal B^{q,d,p}(x-x_0,t; {C_l})$ centered at $x=x_0$ with a very small mass parameter ${C_l}$. In this way we may ensure that $u(x,t)\ge \mathcal B_l(x,t)$ on the parabolic boundary of $E_{t_1}$. Since $\mathcal B_l(x,t)=0$ on $E_{t_1}\cap \{t=0\}$, we only have to check the lateral boundary and this holds if ${C_l}$ is small enough. By the comparison principle on exterior domain for weak or viscosity solutions
$$
u(x,t)\ge \mathcal B_l(x,t) \qquad \mbox{for } \ x\in \ren, \ t=t_1,
$$
since by choosing small enough mass parameter, we can also guarantee that the inequality holds on $B_R(x_0)\times \{ t=t_1 \}$.
The same inequality is then true for $t>t_1$. The conclusion follows.
\qed

The  estimates of the previous theorem are sharp for the explicit solutions by a  direct inspection.

 It is interesting to note the behaviour of the  exponents $\alpha$ and $\alpha/d$ for large values of the parameters since they serve to estimate the size of the solution and the spread rate of the support. Thus, when $q$ is very large, while $p$ and $n$ remain fixed, hence much smaller, we have
$$
d\sim \frac{n-1}{p-1}q, \quad \alpha\sim q^{-1},  \quad \frac{\alpha}{d}\sim \frac{p-1}{n-1}\,q^{-2}\,,
$$
which amounts to  very slow diffusion rates in a fictitious high dimension. On the other hand, when $p, q\to\infty$ with $q=k p$ we get
$$
d\sim (n-1)k+1,  \quad \frac{\alpha}{d}\sim \frac{1}{((n-1)k+2)q}\,.
$$
We see that the fictitious dimension tends to constant, and $\alpha\sim c(n,k)/q$.


\subsection{Sharp asymptotic convergence}
Once the asymptotic size of the solutions is estimated  in the whole space, we proceed with the statement and proof of the precise asymptotic behaviour.

\begin{theorem}
\label{thm:Cauchy-asymp}
Let $u$ be a viscosity solution to the Cauchy problem  (\ref{eq:Cauchy-problem})  and the range condition (\ref{eq:range})
holds. Then, there is  Barenblatt type solution as in (\ref{eq:barenblatt-q-d}) and (\ref{eq:barenblatt.fast}) such that
\begin{align}
\label{as.conv}
\lim_{t\to \infty} t^{\alpha}\sup_{x\in \Rn} \abs{u(x,t)-\B(  x,t;C)}=0.
\end{align}
where $\alpha=d/\lambda$. The constant $C>0 $
depends on the initial data in a non-explicit way.
\end{theorem}

\medskip

\noindent {\sc Proof.} We will give  a proof of this result using ideas from \cite{portilheirov12,portilheirov13}
in the first part, and lap number properties in the final argument.

\medskip

\noindent  Step 1: We will also need a version of Alexandrov's reflection principle. Its proof, which also applies to the equation (\ref{eq:qp}), can be found in \cite[Lemma 9.17, Proposition 14.27]{vazquez07}.

\begin{lemma}[Alexandrov's reflection principle]
\label{lem:alexandrov}
Let $u$ be a viscosity solution to the Cauchy problem (\ref{eq:Cauchy-problem})  with compactly supported initial data $u_{0}\in  C_{0}(B_R(0)), u_{0}\ge 0$. Then for all $t\ge 0$ and all $r>R$ it holds that
\begin{align*}
\min_{\abs x=r} u(x,t) \ge \max_{\abs x=r+2R} u(x,t).
\end{align*}
Moreover, $r \mapsto u(r\theta,t)$, $\theta\in S^{n-1}$, is nonincreasing for $r>R$.
\end{lemma}

The heuristic idea in the proof is to draw a hyperplane $H$ through $(R,0,\ldots,0)$ that divides $\Rn$ into two half spaces $H_+$ and $H_-$, $H^+$ containing $B_R(0)$. Then,  we compare in $H_+$ the solution $u(x,t)$ and $u(\pi(x),t)$ where $\pi(x)$ is the reflection with respect to the hyperplane.
Since the corresponding initial values are ordered and both solutions take the same values on the hyperplane, then the solutions are ordered by the comparison principle. In particular, letting $\nu=(1,0,\ldots,0)$ this gives
\begin{align*}
u(-r\nu)\ge u((2R+r)\nu).
\end{align*}
Then repeating the argument with respect to the other tangent hyperplanes gives the result.

For the nonincreasiness along rays let $H$ be the hyperplane passing perpendicularly trough the middle point of the segment $[r_1 \theta,r_2 \theta]$ for $r_2>r_1>R$.

\medskip

\noindent Step 2: We already have the rough estimates of Theorem \ref{thm:Cauchy-asymp.a}. In order to get sharper approach we need almost radiality for large times in the form used in \cite{portilheirov12, portilheirov13}. We use the previous statement: If  $\spt(u_0)\subset B_{R}(0)$, then we have for all $t\ge 0$ and all $r>R$ that
\[
\begin{split}
\min_{\abs x= r}u(x,t)\ge \max_{\abs x=r+2R}u(x,t).
\end{split}
\]
In order to reduce the gap $2R$, we define the family of rescaled solutions
\begin{equation}\label{scaling.k}
u^{\kappa}(x,t)=\kappa^{\frac{d}{\lambda}}u(\kappa^{\frac1\lambda} x,\kappa t)
\end{equation}
for variable scaling parameter $\kappa>1$, and similarly $\mathcal B_l^{\kappa},\ \mathcal B_u^{\kappa}$, where $\mathcal B_l, \mathcal B_u$ are as in the proof of Theorem \ref{thm:Cauchy-asymp.a}.  This is a  particular case of transformation \eqref{scaling} where the $d$-mass is conserved. These rescaled functions also solve the  same equation (\ref{eq:qp}), and since $\mathcal B_l$ and $\mathcal B_u$ are invariant under this scaling, it follows that for any compact time interval $T/2\le t\le T$  and $q>2$ there is $R_*>0$ such that the support of
$u^{\kappa}(\cdot,t)$ is contained in $ B_{R_*}$
independent of $\kappa$, and that $u^{\kappa}$ are uniformly bounded as well as continuous. In the case $q\le 2$, the solutions no longer have a compact support but uniform boundedness still holds.

Let $\eps>0$, and $\kappa>1$ large enough so that $R_{\kappa}:=\kappa^{-\frac1\lambda}R=\eps$. Once we apply the above inequality for the rescaled solutions, we get for $r>R_{\kappa}$
$$
\min_{\abs x= r}u^{\kappa}(x,t)\ge \max_{\abs x=r+2R_{\kappa}}u^{\kappa}(x,t).
$$
Next we fix $t=1$ and define
\[
\begin{split}
u^{\kappa}_{\min}(r)=\min_{\abs x= r} u^{\kappa}(x,1),\qquad u^{\kappa}_{\max}(r)=\max_{\abs x= r} u^{\kappa}(x,1).
\end{split}
\]
By Lemma \ref{lem:alexandrov}, they are nonincreasing  functions for  $r\ge \eps$, compactly supported in $B_{R_*}(0)$, and satisfy
\[
\begin{split}
 u^{\kappa}_{\max}(r)\ge u^{\kappa}_{\min}(r) \ge u^{\kappa}_{\max}(r+2\eps)
\end{split}
\]
for all $r\ge \eps$. From this, the boundedness, compact support, and uniform regularity for solutions \cite{imbertjs}, it follows that there is a small error in $d$-mass:
\begin{equation}
\label{eq:order}
\begin{split}
\int_0^\infty (u^{\kappa}_{\max}(r)- u^{\kappa}_{\min}(r))r^{d-1}\ud r\le C \eps .
\end{split}
\end{equation}
We have thus arrived at a small asymptotic error in some integral norm  that is propagated in time. If $q<2$, combine the above argument with the tail estimate obtained from \eqref{eq:tail} to obtain (\ref{eq:order}).

\medskip

\noindent Step 3:  Let $u^{\kappa}_{\min}(r,t)$ and $u^{\kappa}_{\max}(r,t)$ be the radial solutions with the initial data given by $u^{\kappa}_{\min}(r)$ and $u^{\kappa}_{\max}(r)$ respectively at $t=1$. Then by the comparison principle for all $t\ge 1$, we have
\begin{equation}
\label{eq:sandwich}
\begin{split}
u^{\kappa}_{\min}(r,t)\le u^{\kappa}(x,t) \le u^{\kappa}_{\max}(r,t)\,,
\end{split}
\end{equation}
where of course $r=\abs x$.

Next we need to find asymptotics for the radial solutions $u^{\kappa}_{\min}(r,t), u^{\kappa}_{\max}(r,t)$. If $d$ were an integer, we  could use the results of \cite{kaminv88}, but since this is not necessarily the case we derive the results for our equation in the next section. Let $\eta>0$. Inspecting the proof of Proposition \ref{prop:Linfty}, as we can always choose $\kappa$ and $t_0$ large enough so that $\kappa=t_0/2$ where $t_0$ is as in the proof of Proposition \ref{prop:Linfty},  we get
\[
\begin{split}
\abs{u^{\kappa}_{\min}(r,2)-\mathcal B^{q,d,p}(r,2,C_{\min})}\le \eta,\\
\abs{u^{\kappa}_{\max}(r,2)-\mathcal B^{q,d,p}(r,2,C_{\max})}\le \eta.
\end{split}
\]
Going back to the original scaling similarly as in the proof of Proposition \ref{prop:Linfty}, we get asymptotic limits for the radial upper and lower bounds in (\ref{eq:sandwich}), and further recalling (\ref{eq:order}) to see that $C_{\max}-C_{\min}$ is small, we obtain the claim.
\qed

\subsection{Convergence of radial solutions via intersection comparison}

For our general equation $d$ is not necessarily an integer, so  the last item of the above proof must be completed. We have to write down the proof of asymptotic convergence for the 1D  equation \eqref{eq.qd} that is satisfied by the radial solutions in the fictitious dimension. This merits a careful consideration and explanation.

We could establish the result about the $d$-mass convergence as in the previous proof. There are many  techniques that have been used to establish the pointwise large time asymptotic behaviour of solutions of parabolic equations. The one we use here works in one dimension, or for radially symmetric solutions in several dimensions, and is based on counting the evolution in time of the ``number of intersections of two solutions'', a rough idea that can be made precise with the names  {\sl intersection number} or {\sl lap number}. These concepts have been investigated in works by Sattinger \cite{Sat69}, Matano \cite{Mat82},  Angenent \cite{An88} and others, and were used by Galaktionov and the second author in a number of papers, cf.\ \cite{GalVazBk}. The idea seems to go back to Sturm, \cite{St36}, so it is also called Sturmian theory, \cite{Gal04}.

First we consider functions $w(x)$ of the real variable $x \in \re$ and look for a count of the number of sign changes by looking for finite sequences of points $x_1<x_2<\cdots<x_{k+1}$  such that $w(x_j)w(x_{j+1})<0$. This is called a sign-change sequence. We define the counter
\begin{equation}
I(w)=\sup\{ k\in \N : \mbox{\rm there exists a sign-change sequence} \ x_1<x_2<\cdots<x_{k+1} \}
\end{equation}
and if there is no sign change we set $I(w)=0$.
Hence, the counter is a nonnegative integer or plus infinity.
In the case of two functions, we denote
\begin{align*}
N(t,u_1,u_2)=I(u_1(\cdot,t)-u_2(\cdot,t)).
\end{align*}
The result we will use is the following improvement of the usual Maximum Principle.

\medskip

Replacing the maximum principle by the elliptic (i.e.\ the one stated with elliptic or Euclidean boundary instead of parabolic boundary) type comparison principle  in the proof of Theorem 4 in \cite{Sat69}, we obtain the Sturmian comparison principle.

 \begin{theorem}[Sturmian Comparison Principle for parabolic equations]
Let $u$ and $u_2$ be two viscosity solutions with possibly different initial data to the Cauchy problem (\ref{eq:Cauchy-problem}).  Then the counter $N(t,u,u_2)$ does not increase in time.
\end{theorem}

This explains the name intersection comparison in the name of the section. We have to specify the equations to which the Sturmian comparison principle applies. The original applications concerned solutions of the classical heat equation in 1D or with radial symmetry in several space dimensions. The application we use here is taken from \cite{vazquez03} where it is applied to the radial solutions of the porous medium equation. It is known to apply for example to $p$-Laplace equations \cite{GenBr07} and reaction-diffusion equations    \cite{QRV02}.

Let us first note  that when the initial counter is zero for the difference of two solutions, then it is zero for all times and the two solutions being compared are ordered at all times. This is a version of the usual comparison principle. The case that interests us is when we consider two solutions, one of them is the solution under investigation, the other one is the Barenblatt solution with the \emph{same  $d$-mass},  and we have  $N(0,u,u_2)=1$, and  then the counter must be 1 or zero for all later times. But the result $N(0,u,u_2)=0$  would imply ordering, and by virtue of the conservation of  $d$-mass this means that the two solutions must be the same at that time, hence the same for all later times, and the asymptotic behaviour is proved.

We have to examine the case where  $N(t,u,u_2)=1$  for all times, and try to make the asymptotic conclusion also in that case. We copy Lemma 17.2 from \cite{vazquez03}. The lemma also holds for our equation. In particular, the fact that the support of any nonnegative $\not\equiv0$ solution spreads with time and occupies any fixed ball follows by comparing with the Barenblatt type solution in the case $q>2$. The case $q\le 2$ can be completed using the ideas of Theorem \ref{thm:Cauchy-asymp.a}.

\begin{lemma}
\label{lem:change-once}
Let $u$ be a radial viscosity solution to the Cauchy problem  (\ref{eq:Cauchy-problem})  and the range condition (\ref{eq:range})
holds. Then there exist time delays $0  <  t_1 < t_2$ such that if
$u_i(r, t) = \mathcal B^{q,d,p}(r, t + t_i;C )$, $i = 1, 2,$ and $C$ is the $d$-mass of $u$, we have
\begin{equation}\label{int.comp.n=1}
N(t, u, u_1) = N(t, u, u_2) = 1.
\end{equation}
Moreover,  we may choose the delays so that for  all large $t$ the list of signs of $u - u_1$ is $- +$, while the list for $u- u_2 $ is $+ -$. Therefore, for
all $r$  small enough we have for all large $ t$
\begin{equation}\label{int.comp}
u_2(r, t) < u(r, t) < u_1(r, t) .
\end{equation}
\end{lemma}

For the conclusion about the convergence of the solutions we need to pay  more attention to formula \eqref{int.comp.n=1} and look for the relative position of the intersections of the three solutions.
Among other things the proof utilizes the conservation of mass. As pointed out at the end of Section \ref{sec.1dim}, the conservation of mass also holds in the singular range as long as the range condition (\ref{eq:range}) holds.


\begin{proposition}
Let $u$ be as in the previous lemma. Then, there is  $C>0$ such that
\begin{equation}
\lim_{t\to \infty} \int_0^\infty \abs{u(r,t)-\B(r,t;C)}r^{d-1}dr=0,
\end{equation}
and such that $u$ and $\B$ have the same $d$-mass.
\end{proposition}

\noindent {\sc Proof.}
Let $u_1$ and $u_2$ be as in the previous lemma so that  $u,u_1,u_2$ have the same $d$-mass.
 For any $t>0$ we denote by $z_1$ the intersection of $u(\cdot,t)$ and $u_1(\cdot,t)$,  by $z_2$ the intersection of $u(\cdot,t)$ and $u_2(\cdot,t)$,  and by $z$ the intersection of $u_1(\cdot,t)$ and $u_2(\cdot,t)$.   By the formulas of the Barenblatt solutions we see that $N(t,u_1,u_2)=1$  and also observe that $z$ is known a priori. There are three cases to consider.

\noindent{\bf First case:}  Suppose that $z_1>z$. Looking at the graphs of the functions we see that $u$ must have gone under $u_2(r,t)$ before reaching $r=z$, and then $z_2\le z$. This means that $u$ lies between the two solutions for $0\le r\le z_2$ and $r\ge z_1$, and moreover,
$$
u(r,t)\le \min\{u_1(r,t),  u_2(r,t) \} \quad \text{ for } \ z_2\le r\le z_1.
$$

\noindent{\bf Second case:}  Suppose that $z_2<z=z_1$. This is similar to the first case.

\noindent{\bf Third case:}  Finally, $z_1\le z\le z_2$.   In this case we have
\begin{align*}
u_2(r,t)&\le u(r,t)\le u_1(r,t) \text{ for }0<r\le z_1\\
u_1(r,t)&\le u(r,t)\le u_2(r,t) \text{ for }z_2\le r<\infty.
\end{align*}

Moreover,
$$
u(r,t)\ge \max\{u_1(r,t),u_2(r,t)\}\text{ for }z_1\le r\le z_2.
$$
\sloppy
\noindent {\bf Asymptotic mass analysis.} We recall that  $u(r,t), u_1(r,t),u_2(r,t)$, as well as $\B( r,t;C)$, have the same $d$-mass for all times.  We want to prove that
$$
\lim_{t\to \infty} \int_0^\infty \abs{u(r,t)-\B( r,t;C)}r^{d-1}dr=0.
$$
Let us assume that we are in the third case above, other cases being similar. We estimate the integral in the three regions $I_1=[0,z_1]$,
$I_1=[z_1,z_2]$, and $I_3=[z_2,\infty]$. In $I_1$ and $I_3$ we have $u$ sandwiched between $u_1$ and $u_2$ so that
$$
\abs{u(r,t)-u_1(r,t)}\le  \abs{u_2(r,t)-u_1(r,t)},
$$
and moreover it can be easily seen that the $d$-mass of the  right hand side goes to zero. Assume that the integrals on both  intervals  are  less than or equal to $\ve /4$ for $t\ge t_\ve$. We now estimate the integral in the middle interval as follows:
\begin{align*}
 \int_{I_2} &\abs{u(r,t)-u_1( r,t)}r^{d-1}dr= \int_{I_2} ( u(r,t)-u_1( r,t))\,r^{d-1}dr\\
&=\int_0^{\infty} (u(r,t)-u_1( r,t))\,r^{d-1}dr- \int_{I_1\cup I_3} (u(r,t)-u_1( r,t))\,r^{d-1}dr\le 0+ 2\ve/4.
\end{align*}
Here, we used the order in the interval $I_2$ and the equality of the total $d$-mass. The claim follows.
 \qed

\begin{proposition}
\label{prop:Linfty}
Let $u$ be as in the previous lemma. Then, there is $C>0$ such that
\begin{equation}
\lim_{t\to \infty} t^{d/\lambda}\abs{u(r,t)-\B(r,t;C)}=0.
\end{equation}
\end{proposition}

\noindent {\sc Proof.} Let $u_1$, and $u_2$ be as in the previous lemma, and let  $u,u_1,u_2$ have the same $d$-mass. By scaling we may assume that $t=1$ and that the $d$-mass of the difference $\abs{u(r,1)-\B(r,1;C)}$ is small. To be more precise, we first choose $t_0$ large enough so that the $d$-mass of the absolute value of the difference is small. Then we consider the rescaled function $u(r,t):=u^{\kappa}(r,t)$ in \eqref{scaling.k} with $\kappa=t_0$. Also observe that $d$-mass of the difference remains small in this rescaling.

Since we know that the solutions are uniformly bounded in $L^\infty(\re)$ by a comparison argument similar to those in Theorem \ref{thm:Cauchy-asymp.a},
 the $C^1$ regularity result, Theorem 1.1 in \cite{imbertjs}, implies that the derivative of  $u(r,1)-\B(r,1;C)$ is uniformly bounded. But then a simple argument says that $u(r,1)-\B(r,1;C)$  is small in $L^\infty(\re)$.

By returning to the original scaling, we obtain the power factor $ t^{d/\lambda}$. \qed


\section{Asymptotic behaviour in bounded domains}
\label{sec:bounded-dom}

In this section, we consider the equation defined in a bounded domain with zero Dirichlet boundary conditions. We consider  the problem
\begin{equation}
\label{eq:dir-prob}
\begin{split}
\begin{cases}
u_t=\abs{Du}^{q-2}(\Delta u+(p-2)\Delta_{\infty}^{N} u),& \text{ in } \Om_{\infty},\\
u(x,0)=u_0(x),& x\in \Omega,\\
u(x,t)=0,& x\in \partial \Omega,
\end{cases}
\end{split}
\end{equation}
where  $\Om_{\infty}=\Omega\times (0,\infty)$ where  $\Omega$ is a bounded domain of $\ren$ with smooth enough boundary.
We will always take  $q>2$ and $p> 1$, and $u_0\ge 0$, $u_0\not\equiv 0$ bounded $C^2(\ol \Om)$-function.

The study of the behaviour of nonnegative solutions defined in a general bounded domain  is based on three ingredients: the existence of uniform bounds, the study of the special solution on a ball, and a monotonicity condition that  holds for nonnegative solutions. We will be specially interested in the long-time behaviour.

%

\subsection{Properties of solutions defined in a general bounded domain}

We will follow the outline of the proof of \cite{Vasbdd} for the porous medium equation with a number of changes needed in our framework, so the proof is differently organized. We note that some arguments are only sketched whenever  they can be easily adapted.

\medskip

\noindent $\bullet$  {\sl Monotonicity} is a key property in our proofs of long-time behaviour. It takes the form of a derivative bound that is in fact valid for all semigroups generated by a homogeneous operator acting on a Banach space or a convex subset thereof. We state it as follows

\begin{proposition} The nonnegative viscosity solutions of  (\ref{eq:dir-prob})  satisfy
\begin{equation}\label{monot.c}
\partial_t u\ge -\frac{u}{(q-2)t}\,.
\end{equation}
in the sense of distributions. Recall that  $q>2$.
\end{proposition}

\noindent {\sl Sketch of proof.}  The original proof is done in   \cite{BC1981}, see also  Lemma 8.1 in \cite{vazquez07}. We briefly repeat the ideas: we denote by $S_tu_0 $ the solution with initial data $u_0$. Now we recall that the rescaling implies that if $u(x,t)$ is a solution then $u_k(x,t)=k u(x,k^{q-2}t)$ is again a solution for all $k>0$. This is equivalent to writing  $S_t(k u_0)= k S_{k^{q-2}t}(u_0)$, or putting $\lambda=k^{q-2}$:
$$
S_t(\lambda^{\frac1{q-2}} u_0)= \lambda^{\frac1{q-2}}  S_{\lambda t}(u_0).
$$
 Then,
$$
\begin{array}{c}
S_{\lambda t}(u_0)-S_{t}(u_0)= \lambda^{-\frac1{q-2}}S_t(\lambda^{\frac1{q-2}} u_0)-S_{t}(u_0)\\
=(\lambda^{-\frac1{q-2}}-1)S_t(\lambda^{\frac1{q-2}} u_0)+ ( S_t(\lambda^{\frac1{q-2}} u_0)-S_{t}(u_0)).
\end{array}
$$
By the maximum principle the last  summand is positive whenever $\lambda>1$. Now observe that for $\lambda=1+\ve>1$ we have
for $\ve >0$ small
$$
0\ge \lambda^{-\frac1{q-2}}-1= -(1-(1+\ve)^{-\frac1{q-2}})\sim -\frac1{q-2}\ve
$$
while $\lambda t- t\sim  t \ve$. Taking the incremental quotient, $(S_{\lambda t}(u_0)-S_{t}(u_0))/(\lambda t-t)$,  and passing to the limit $\ve\to 0$ the result holds.

 Note that the proof works for the viscosity solutions: observe that multiplying by a test function above, integrating, and moving the difference quotient on the test function, we get the statement in the sense of distributions.
\qed

\medskip

\noindent $\bullet$ Another main ingredient will be {\sl rescaling} and the study of a {\sl new differential equation}.  We perform the following change of variables
\begin{equation}\label{eq:resc}
u(x,t)=t^{-\frac1{q-2}}v(x,\tau), \qquad \tau=\log(t)\,.
\end{equation}
The equation for $v$ is
\begin{align}
\label{eq:resc-eq}
\partial_\tau v= \abs{Dv}^{q-2}(\Delta v+(p-2)\Delta_{\infty}^{N} v)+\frac1{q-2} v.
\end{align}
It follows that $v(x,\tau)\ge 0$ is defined for all $\tau>-\infty$.   Moreover, we will see below that $v$ is uniformly bounded in $x$ and $\tau$ as a consequence of the construction of the next subsection.

\subsection{Asymptotic behaviour in a ball }

We consider now in detail the problem posed in a  ball $\Omega=B_R(0)$.  By scaling we may take $R=1$ without loss of generality. Indeed, if $u(x,t)$ is a solution with space domain $B_1(0)$, then \
$u_R(x,t)=Au(x/R, t)$ is another solution defined in $B_R(0)$ if and only if $A^{q-2}=R^{q}$ .

 \medskip

\noindent {\bf Separable solutions.} We study radial separable solutions with zero lateral boundary values having a form
\begin{equation}
\label{eq:separable}
\begin{split}
U(x,t)=t^{-\frac1{q-2}}V(\abs{x})
\end{split}
\end{equation}
sometimes called the friendly giant, or FG-type solution. This is a particular case  of (\ref{eq:resc}) where the second factor is time-independent.

\begin{theorem}\label{thm.6.2}
The separable weak/viscosity solution (\ref{eq:separable}) defined in $B_1(0)$ exists and is unique in the class of positive profiles $V$. 
\end{theorem}

\noindent {\sl Proof.}  { Step 1: Preparation}. Using the equation, recalling (\ref{eq:radial-p-Laplace}) and (\ref{eq:1-divergence-form}), we are led to consider  the stationary problem for $V$
\begin{equation}
\label{eq:elliptic}
\begin{split}
0=\frac{V(y)}{q-2}\abs y^{d-1}+\frac{p-1}{q-1}(\abs{V'(y)}^{q-2}V'(y) \abs y^{d-1})',
\end{split}
\end{equation}
for $y\in (-1,1)$. Using  symmetry we will  look for a bounded solution  $V$ of the equation only in $0<y<1$  with the conditions $V\ge 0, V'\le 0$, and end point values
\[
\lim_{y\to 0+} \abs{V'(y)}^{q-1} y^{d-1}=0,\quad   V(1)=0.
\]
  In order to prove the existence, we write  for $r\in [0,1]$
\[
\begin{split}
-\int_{0}^r \frac{V(y)}{q-2} y^{d-1} \ud y & =\frac{p-1}{q-1}\abs{V'(r)}^{q-2}V'(r)  r^{d-1}
\end{split}
\]
i.e.
\begin{align}
\label{eq:grad-reg}
-\Bigg(\frac{q-1}{p-1}\int_{0}^r \frac{V(y)}{q-2}\Big(\frac{y}{r}\Big)^{d-1} \ud y\Bigg)^{\frac1{q-1}}=V'(r),
\end{align}
and further  for $R\in [0,1]$
\begin{equation}
\label{eq:int-eq}
\begin{split}
\int^1_R\Bigg(\frac{q-1}{p-1}\int_{0}^r \frac{V(y)}{q-2}\Big(\frac{y}{r}\Big)^{d-1} \ud y\Bigg)^{\frac1{q-1}}\ud r=V(R).
\end{split}
\end{equation}
From here it follows that if we have a solution $V\in C([0,1])$  for (\ref{eq:int-eq}), then $V\in C^1([0,1])$ and $\lim_{r\to 0+}V'(r)=0$. Then reflecting $V$  evenly to define $V$ over the whole interval $(-1,1)$, the function $U$  given by \eqref{eq:separable} is a weak solution according to Definition \ref{def:1D-weak}. Theorem \ref{thm:equiv} shows that we have also constructed a viscosity solution to our problem.

\medskip

\noindent  {Step 2: Existence.} It remains to verify that (\ref{eq:int-eq}) has a solution. Let
\[
\begin{split}
\mathcal C&=\{ u\in C([0,1]) \,:\,u(1)=0,u\ge 0\},\\
\mathcal K&=\{ u\in \mathcal C \,:\, \norm{u}_{L^\infty(0,1)}\le M,\norm{u}_{L^\infty(0,\half)}\ge m,\ \abs{u(x)-u(y)}\le M\abs{x-y},\  x,y\in [0,1]\}
\end{split}
\]
where the $q$ dependent constants $M:=M(q,d)$  and $m:=m(q,d)$ will be determined in the course of the proof. The left hand side of (\ref{eq:int-eq}) defines an operator $T$.
 Since $\mathcal K$ is convex, closed and compact (by Arzel\`a-Ascoli), we can use Schauder's fixed point theorem to find a nontrivial solution if $T:\mathcal K\to \mathcal K$ is a continuous operator.
To verify this, observe that $TV(1)=0$ and that for $u\in \mathcal K$ it holds that
\[
\begin{split}
\abs{TV(R)}&=\abs{\int^1_R\Bigg(\frac{q-1}{p-1}\int_{0}^r \frac{V(y)}{q-2}\Big(\frac{y}{r}\Big)^{d-1} \ud y\Bigg)^{\frac1{q-1}}\ud r}\\
&\le M^{\frac1{q-1}}\abs{\int^1_R\Bigg(\frac{q-1}{p-1}\int_{0}^r \frac{1}{q-2}\Big(\frac{y}{r}\Big)^{d-1} \ud y\Bigg)^{\frac1{q-1}}\ud r}\\
&= M^{\frac1{q-1}}\abs{\int^1_R\Bigg(\frac{r(q-1)}{d(q-2)(p-1)}\Bigg)^{\frac1{q-1}}\ud r}\\
&\le  \Bigg(\frac{M(q-1)}{d(q-2)(p-1)}\Bigg)^{\frac1{q-1}} (1-R)\le M
\end{split}
\]
since $r\le 1$ and the last inequality holds for large enough $M=M(q,d)$. This is because $q>2$, $M^{\frac1{q-1}} =M M^{\frac{2-q}{q-1}}$ and  $M^{\frac{2-q}{q-1}}$ can be made small by choosing large enough $M$. Similarly, whenever $R\in[0,\half]$, it holds that
 \[
\begin{split}
\abs{TV(R)}&\ge m^{\frac1{q-1}}\abs{\int^1_{\half}\Bigg(\frac{r(q-1)}{d(q-2)(p-1)}\Bigg)^{\frac1{q-1}}\ud r}\ge m
\end{split}
\]
for small enough $m=m(q,d)>0$, and in particular we get a nontrivial solution.
 Moreover, assume without loss of generality that $R_1\le R_2$ and observe  by a similar computation as above that whenever $V$ is Lipschitz with a constant $M$ we have
 \[
\begin{split}
\abs{TV(R_2)-TV(R_1)}=\abs{\int^{R_2}_{R_1}\Bigg(\frac{q-1}{p-1}\int_{0}^r \frac{V(y)}{q-2}\Big(\frac{y}{r}\Big)^{d-1} \ud y\Bigg)^{\frac1{q-1}}\ud r}\le M\abs{R_2-R_1}.
\end{split}
\]
Finally,
\[
\begin{split}
&\norm{\int^1_R\Bigg(\frac{q-1}{p-1}\int_{0}^r \frac{V_1(y)}{q-2}\Big(\frac{y}{r}\Big)^{d-1} \ud y\Bigg)^{\frac1{q-1}}-\Bigg(\frac{q-1}{p-1}\int_{0}^r \frac{V_2(y)}{q-2}\Big(\frac{y}{r}\Big)^{d-1} \ud y\Bigg)^{\frac1{q-1}}\ud r}_{L^{\infty}(0,1)}\\
&\le \norm{\int^1_R\Big(\frac{q-1}{p-1}\int_{0}^r \frac{\abs{V_1(y)-V_2(y)}}{q-2}\Big(\frac{y}{r}\Big)^{d-1} \ud y\Big)^{\frac1{q-1}}\ud r}_{L^{\infty}(0,1)}\le C\norm{V_1-V_2}_{L^{\infty}(0,1)}^{\frac1{q-1}},
\end{split}
\]
and thus $T$ is continuous.

\medskip

\noindent {Step 3: Uniqueness.} Note first that by scaling we can construct a separable solution in a domain $B_R(0)$ for any $R>0$ by means of the formula
\begin{equation}\label{VR.eq}
U_R(r,t)=t^{-\frac1{q-2}}V_R(r)\,,\qquad V_R(r)=R^{q/(q-2)}V(r/R), \qquad 0<r<R.
\end{equation}
Suppose now that $U$ and $\widetilde U$ are two solutions with profiles $V(r)$ and
$\widetilde V(r)$ respectively. First we construct the solution $U_{1+\eps}$ of the separable form  for domain $B_{1+\eps}(0)\times (0,\infty)$ by using on $U$ the above rescaling from radius $1$ to $1+\eps$. Taking a small enough $t_0>0$ it is easy to see  that $\widetilde U (r,1)\le U_{1+\eps}(r,t_0)$ in $B_1(0)$. By the parabolic comparison principle  we get for every $t>0$ and $0<r<1$
$$
\widetilde U (r,1+t)\le U_{1+\eps}(r,t_0+t)\,,
$$
i.e.
$$
\widetilde V (r)\le \left(\frac{1+t}{t_0+t}\right)^{1/(q-2)}V_{1+\eps}(r)\,.
$$
We pass to the limit $t\to\infty$ to get $\widetilde V (r)\le V_{1+\eps}(r)$. Now let $\eps\to 0$ and use the scaling law (\ref{VR.eq}) to get $\widetilde V (r)\le V(r)$. The other inequality is obtained in the same way, hence $\widetilde V (r)= V(r)$. \qed

The next result follows from the above considerations by using  the equivalence result of Remark~\ref{rem:extensions} (ii), or the argument of Theorem 7.2 in \cite{portilheirov12} as in Step 6 of the next section.

\begin{corollary}\label{coroll.statsol} The  profile $V>0$ is the unique  nonnegative viscosity solution of the stationary problem
\begin{equation}
\abs{DV}^{q-2}(\Delta V+(p-2)\Delta_{\infty}^{N} V)+\frac1{q-2} V=0,\quad \text{ in } \  B_1(0),
\end{equation}
with zero Dirichlet boundary conditions.We have $V\in C^{1,\beta}(-1,1),\ \beta\in (0,1)$.
\end{corollary}

\noindent {\bf Remark on the fast case.}
If $q\le 2$, then we obtain no friendly giant type solution but
\[
\begin{split}
u(x,t)=
\begin{cases}
 (t_*-t)^{-\frac1{q-2}}V(x)& 0\le t< t_*\\
0 &  t\ge t_*
\end{cases}
\end{split}
\]
is a solution with a suitable $V$ as shown in \cite{ohnumas97}, Section 5, in 1-dimensional case. This shows that the threshold in the bounded domain case for extinction in finite time is $q=2$.

\medskip

\noindent $\bullet$ Using the above separable solution $U$, we obtain the decay rate of general radial solutions as $u(x,t)=O(t^{-\frac1{q-2}})$, and much more: we also obtain the precise asymptotic behaviour of radial solutions.

\begin{theorem}
\label{thm:bounded-dom}
Suppose that $u_0\in C^2(\ol B_1), u_0(x)\ge 0$, $u_0\not\equiv 0$ for $x\in B_1(0)$, $u_0$ radial and $u_0(x)=0$ on $x\in \partial B_1(0)$. Then, the solution to (\ref{eq:dir-prob}) satisfies
\begin{equation}\label{conv.rad}
\lim_{t\to \infty}t^{\frac1{q-2}}u(x,t)=V(\abs x)\,,
\end{equation}
uniformly in $x\in B_1(0)$, and also
\begin{equation}
u(x,t)\le U(x,t)=t^{-\frac1{q-2}}V(\abs{x})
\end{equation}
for every $x\in B_1(0)$, $t>0$. The positive bounded function $V$ is the stationary solution  of Corollary \ref{coroll.statsol}.
\end{theorem}

\noindent {\sc Proof.} {Step 1: Universal Boundedness.}  Arguing like in the  comparison argument used in the above uniqueness proof, we see  that for every $\eps>0$ there is $t_0>0$ such that $u(x,0)\le t_0^{-1/(q-2)}V_{1+\eps}(\abs x)$, so that the comparison argument implies that
$$
u(x,t)\le (t+t_0)^{-1/(q-2)}V_{1+\eps}(\abs x)\le t^{-1/(q-2)}V_{1+\eps}(\abs x)\,.
$$
This is a uniform bound on all solutions of the problem.

\medskip

 \noindent {Step 2: Rescaling and new equation}. As indicated above in (\ref{eq:resc}), we perform the  change of variables
$u(x,t)=t^{-\frac1{q-2}}v(x,\tau),$ with $\tau=\log(t),$ we get the equation for $v$
\begin{align}
\label{eq:resc-eq}
\partial_\tau v= \abs{Dv}^{q-2}(\Delta v+(p-2)\Delta_{\infty}^{N} v)+\frac1{q-2} v.
\end{align}
It follows that $v(x,\tau)\ge 0$ is defined for all $\tau>-\infty$. It is uniformly bounded in $x$ and $\tau$ by virtue of the universal estimate for $u$ we have just proved.

\noindent {Step 3: The limit exists and is positive as well as bounded.} The monotonicity condition \eqref{monot.c} is equivalent to $\partial_\tau v\ge 0$  in distributional sense. This and boundedness imply that there exists a limit
\begin{equation}\label{eq:monot-limit}
\lim_{\tau\to \infty}v(x,\tau)=W(x)\le V(|x|)
\end{equation}
at least in a pointwise sense in $B_1(0)$, since the limit is independent of $\eps$. We have $W(x)\ge 0$. In fact $W(x)$ is strictly positive in $\Omega$
by a comparison argument applied to $v(x,\tau)$ by comparing with small Barenblatt solutions used as subsolutions. The comparison first proves that a point of positivity of the solution of $v(x,\tau_0)$ stays positive for $v(x,\tau)$ with $\tau>\tau_0$. But it also proves that the positivity set of $v(x,\tau)$ expands with time to cover all points of $\Omega$, a connected set. We conclude that $W(x)$ must be positive everywhere.

\noindent {Step 4: Identification of the limit.} Now we perform the comparison of Step (i) in the other direction. Given our solution and any $\eps>0$ we can find a large $t_1$ so that the positivity set of $u$ covers $B_{1-\eps/2}(0)$. Then we choose $t_2$ large enough so that
$$
u(x,t_1)\ge t_2^{-1/(q-2)}V_{1-\eps}(r) \quad \mbox{ on} \ B_{1-\eps}(0).
 $$
Taking these as initial functions, it easily follows by comparison of  viscosity solutions in $B_{1-\eps}(0)\times (t_1,\infty)$ that
$$
u(x,t+t_1)\ge (t+t_2)^{-1/(q-2)}V_{1-\eps}(r).
$$
Note that the ordering also holds on $\partial B_{1-\eps}(0)$. Passing  in this inequality to the limit $t\to\infty$ we get
$$
\lim_{t\to \infty}t^{\frac1{q-2}}u(x,t)\ge V_{1-\eps}(r).
$$
Let now $\eps\to0$ and we get the convergence result \eqref{conv.rad} of the theorem. The regularity is a consequence of the existence construction, see (\ref{eq:grad-reg}), and thus  the proof is complete.\qed

\medskip

The result also implies that  $U$ is the minimal universal upper bound of the class of solutions.

\subsection{Asymptotic behaviour in a general bounded domain}

The study of the asymptotic behaviour in a general bounded domain, though more difficult,  is based on similar ideas.

\begin{theorem}
\label{thm:bounded-dom} (i) Let $\Omega$ be a bounded domain in $\ren$ with $C^2$ boundary and $q>2$.
Let $u$ be a  viscosity solution to \eqref{eq:dir-prob} posed in $\Om_{\infty}=\Omega\times(0,\infty)$ with zero Dirichlet lateral boundary conditions, and initial data $u_0\in C^2(\ol \Omega)$, $u_0(x)>0$ for $x\in \Omega$, and $u_0(x)=0$ on $x\in \partial \Omega$. Then,  we have
\begin{equation}\label{conv.monot}
\lim_{t\to \infty}t^{\frac1{q-2}}u(x,t)=W(x)
\end{equation}
uniformly in $x$, where $W$ is strictly positive and bounded.  Moreover, $u(x,t)\le U(x,t):=t^{-\frac1{q-2}}W(x)$ for every $x\in\ren$, $t>0$.

(ii) If $\Om_1\Subset \Om_2$ and $W_1, W_2$ are the respective limits, then $W_1\le W_2$.

(iii) If moreover $\Omega$ is starshaped then the  profile $W\in C^{1,\beta}(\Omega)$ is the unique positive viscosity solution of the stationary problem
\begin{equation}
\abs{DW}^{q-2}(\Delta W+(p-2)\Delta_{\infty}^{N} W)+\frac1{q-2} W=0,\quad \text{ in } \ \Omega
\end{equation}
with zero Dirichlet boundary conditions. The expression
\begin{equation}
U(x,t)=t^{-\frac1{q-2}}W(x)
\end{equation}
is a particular viscosity solution  to \eqref{eq:dir-prob} posed in $\Om_{\infty} $ with zero Dirichlet lateral boundary conditions. Note that $U$ takes infinite initial data.
\end{theorem}

\noindent {\sl Proof.}  The Steps 1-4 below work for general bounded domains.

\noindent {Step 1: Universal Boundedness.} This  follows from comparison  with the friendly giant type solution of the problem posed in a larger ball, $\Omega\subset B_R(x_0)$. Let us call the $x$ dependent part $V_R(|x-x_0|)$. An easy comparison shows that for any solution $u$  of our problem we have
$$
u(x,t)\le t^{-\frac1{q-2}}V_R(\abs{x-x_0})\le C\,t^{-\frac1{q-2}}.
$$
This is indeed a universal estimate.

\noindent {Step 2: Rescaling and new equation}.  This step is identical to the previous proof.

\noindent {Step 3: The limit (\ref{conv.monot}) exists and is  positive as well as bounded.} From the monotonicity we get the existence of the  limit $W$ and it is a nonnegative, bounded function.  The fact that $W$ is strictly positive everywhere in $\Omega$ follows by a comparison argument applied to $v(x,\tau)$. It works exactly as in the radial case  by comparing with small Barenblatt solutions used as subsolutions during the period of time until the support
reaches the boundary (this is called expansion of the positivity set). We have proven (i).

\noindent {Step 4: The limit is monotone with respect to the domain.} Suppose that $\Omega_1\subset \Omega_2$ with a positive distance from $\partial \Omega_1$ to $\partial \Omega_2$. Let $u_1$ and $u_2$ be solutions as above defined respectively in $\Omega_1$ and $\Omega_2$. Repeating the comparison argument done in the radial case we can get
times $t_0, t_1>0$ such that
$$
u_1(x,t+t_1)\le u_2(x,t) \quad \mbox{ if} \quad t\ge t_0.
$$
Then the respective limits satisfy $W_1(x)\le W_2(x)$. In particular, in order to get  a lower bound for  $W_2$, limit of a  given solution,  $W_1$ may be chosen as the unique limit of the radial case when we take as $\Omega_1$  a  ball strictly contained in $\Omega.$

\medskip

\noindent Step 5: We now introduce the extra condition on the  domain in order to improve the results. Moving the origin of coordinates to the point that serves as basis of starshapedness, we may define the domains \ $\Omega_\lambda=\{ \lambda x: x\in\Omega\}$ \ for all $\lambda>0$. It follows that for $\eps>0$
$$
\Omega_{1-\eps}\subset \Omega \subset\Omega_{1+\eps}.
$$
We also define the rescalings of the solution much as in the radial case. We take
$u_{1-\eps}(x,t)=Au(x/(1-\eps), t)$, $A^{q-2}=(1-\eps)^{q}$, to get  another solution defined in $\Omega_{1-\eps}$. Given two solutions $u$,  $\widehat u$ and any $\eps>0$ we can find a large $t_1$ so that the positivity set of $u$ covers $\Om_{1-\eps/2}(0)$. Then we choose  $t_2$ large enough so that
$$
u(x,t_1)\ge \widehat  u_{1-\eps}(x,t_1+t_2) \quad \mbox{ on} \ \Om_{1-\eps}(0),
 $$
 since $u$ is positive and continuous in the closure of $\Omega_{1-\eps}$ and $\widehat  u_{1-\eps}$ goes to zero as $t\to\infty$. Taking these  as initial functions, it easily follows by comparison of  viscosity solutions in $\Omega_{1-\eps}\times (t_1,\infty)$ that
$$
u(x,t+t_1)\ge \widehat  u_{1-\eps}(x,t+t_1+t_2)\,,
$$
since the ordering on $\partial \Omega_{1-\eps}$ also holds. Using this inequality and passing to the limit $t\to\infty$ we easily get
$$
W(x)\ge \widehat W_{1-\eps}(x).
$$
Let now $\eps\to0$ we get comparison of the limit profile. Reversing the roles we get uniqueness of the limit.

\noindent {Step 6: \sl  $W$ is a stationary viscosity solution.}  This  is obtained using stability principle for viscosity solutions in the uniform convergence, see Theorem 7.3 in \cite{portilheirov12}. The Lipschitz regularity  for the equation (\ref{eq:resc-eq}) used in Theorem 7.3 of \cite{portilheirov12} follows from similar barrier arguments as in Section 4 of \cite{juutinenk06}, see in particular Corollary 4.3 there, and also Section 6 in \cite{portilheirov12}.

Once we have a bounded solution of the elliptic equation, regularity theory, cf.\  \cite{attouchipr17, attouchir}, means that it will be $C^{1,\beta}$  up to the boundary. \qed

\medskip

\noindent {\bf Remarks.} (1) We are not able to prove the uniqueness of positive solutions of the elliptic problem in general domains. That would imply an asymptotic result as complete as in the stated cases.

(2)  Uniqueness of the positive limit profile is known in particular cases, like the standard $p$-Laplacian (case $q=p$), see \cite{diazs87, anane87}  and for example \cite{bellonik02}.

\section{A priori estimates}\label{apriori}

The next lemma is a counterpart of Proposition III.3.1 in \cite{dibenedetto93} or Proposition 4.7 in \cite{urbano08}. However, as we already have the regularity estimates from \cite{imbertjs}, the proof is simpler. We use the lemma in the proof of Harnack's inequality in Theorem \ref{thm:harnack}.  In the lemma and the following corollary, all the reference points for cylinders are the same, and thus we drop them.
\begin{lemma}
\label{lem:key-osc}
Let $u$ be a viscosity solution to (\ref{eq:qp}) in $Q_{R,R^q}$.
For $\gamma \in (0,1)$, there is $C>1$  that can be determined a priori only depending on $n,p,q, \gamma$ such that the following holds.
Suppose that we are given $\om_0>1$ such that for $a_0=(1/\om_0)^{q-2}$
\[
\begin{split}
\osc_{Q_{R ,a_0 R^q}} u\le \om_0.
\end{split}
\]
Define the sequences
 \[
\begin{split}
R_i&=C^{-i}R,\ R_0=R,\\
\om_{i}&= \gamma \om_{i-1},
\end{split}
\]
where $i=1,2,\ldots$. Then for $Q_{R_i,a_i R_i^q}$ where $a_i=(1/\om_i)^{q-2}$
\[
\begin{split}
Q_{R_{i+1},a_{i+1} R_{i+1}^q}\subset Q_{R_i,a_i R_i^q},\quad \osc_{Q_{R_i,a_i R_i^q}}u \le \om_i.
\end{split}
\]
\end{lemma}
\noindent {\sc Proof.}
First observe that $Q_{R_{i+1},a_{i+1} R_{i+1}^q}\subset Q_{R_i,a_i R_i^q}\subset\ldots \subset Q_{R,R^q}$ holds as long as $C$ and $\gamma$ will satisfy $C^q\gamma^{q-2}\ge1$.

The case $i=0$ holds by the assumption.
 Suppose then that the claim holds for $i=k$ i.e.
\[
\begin{split}
\osc_{Q_{R_k,a_kR_k^q}} u\le \om_k.
\end{split}
\]
By setting
\[
\begin{split}
u_k(x,t):=\frac{u(R_kx,a_kR_k^qt)-\inf_{Q_{R_k,a_kR_k^q}}u}{\om_k}
\end{split}
\]
it holds that $\sup_{Q_{1,1}} u_k\le 1$ by the induction assumption.  Then by \cite[Lemma 2.3, Lemma 3.1]{imbertjs}
\[
\begin{split}
\osc_{Q_{R_{k+1},a_{k+1}R_{k+1}^q}} u/\om_k=\osc_{Q_{C^{-1},\gamma^{-(q-2)}C^{-q} }}u_k\le \tilde C(C^{-1}+(\gamma^{-(q-2)}C^{-q})^{1/2})
\end{split}
\]
where $\tilde C=\tilde C(n,p,q)$ is the constant in the regularity estimates  cited above. Choosing  $C>\max\{ 2\tilde C/\gamma , (2\tilde C)^{2/q}\gamma^{-1},\gamma^{-1}\}$ we get
\[
\begin{split}
\osc_{Q_{R_{k+1},a_{k+1}R_{k+1}^q}} u\le \tilde C(\frac{\gamma}{2\tilde C}+\frac{\gamma}{2\tilde C})\om_k=\gamma\om_k=\om_{k+1}.
\end{split}
\]
and $C^{q}\gamma^{q-2}> 1$.
\qed

The standard iteration argument then implies the following corollary.
\begin{corollary}
\label{cor:key-osc}
Let $u$, $\om_0, a_0$, and $R$ be  as in the previous lemma.
Then there exist constants $\hat C=\hat C(n,p,q)>1$ and $\a=\a(n,p,q)\in(0,1)$
 such that for
$
0<r\le R$
it holds that
\[
\begin{split}
\osc_{Q_{r,a_0r^q}} u\le \hat C\om_0\Big(\frac{r}{R}\Big)^{\a}.
\end{split}
\]
\end{corollary}
\noindent {\sc Proof.}
Let $\gamma, C, a_k,R_k$ be as in the previous lemma.
Choose an integer $k$ such that
\[
\begin{split}
C^{-(k+1)}R\le r<C^{-k}R=:R_k.
\end{split}
\]
By this and the recursive definition of $\om_k$ it follows that
\[
\begin{split}
 \om_k=\gamma^{k}\om_0=\gamma^{-1}\gamma^{k+1}\om_0\le\gamma^{-1} \gamma^{-\frac{\log(r/R)}{\log(C)}} \om_0 \le \gamma^{-1}\Big(\frac{r}{R} \Big)^{-\frac{\log(\gamma)}{\log(C)}}\om_0.
\end{split}
\]
By using Lemma \ref{lem:key-osc}, we get
\[
\begin{split}
\osc_{Q_{r,a_0 r^q}} u\le \osc_{Q_{r,a_kr^q}} u\le  \osc_{Q_{R_k,a_kR_k^q}} u \le \om_k.
\end{split}
\]
Setting $\hat C=\gamma^{-1}$ and observing that
$$
\a:=-\frac{\log(\gamma)}{\log(C)}\in (0,1)
$$
since $C>\gamma^{-1}$, the result follows from the previous estimates.
\qed

Next we demonstrate the use of radial solutions and prove Harnack's inequality. The standard proof utilizes the oscillation estimate and expansion of positivity using a radial comparison function.  Thus having the above results at our disposal, the proof is the same as that of Theorem 2.1 on p.\ 157 in \cite{dibenedetto93}; proof is in Chapter VI, Section 4.
\begin{theorem}
\label{thm:harnack}
Let $u\ge 0$ be a viscosity solution to (\ref{eq:qp}) in $Q_{1,1}$ and the range condition (\ref{eq:range}) holds. Fix $(x_0,t_0)\in Q_{1,1}$ and suppose that $u(x_0,t_0)>0$. Then there exist $\mu=\mu(n,p,q)$ and $C=C(n,p,q)$ such that
\[
\begin{split}
u(x_0,t_0)\le \mu \inf_{B_r(x_0)}u(\cdot, t_0+\theta)
\end{split}
\]
where
\[
\begin{split}
\theta=\frac{Cr^{q}}{u(x_0,t_0)^{q-2}},
\end{split}
\]
whenever $B_{4r}(x_0)\times (t_0-4\theta, t_0+4\theta)\subset Q_{1,1}$.
\end{theorem}
\noindent {\sc Proof.} Let first $q>2$.
We consider the rescaled function
\[
\begin{split}
v(x,t)=\frac{1}{u(x_0,t_0)}u\Big(x_0+rx, t_0+\frac{t r^q}{u(x_0,t_0)^{q-2}}\Big)
\end{split}
\]
which is a solution to
\[
\begin{split}
\begin{cases}
v_t=\abs{Dv}^{q-2}(\Delta v+(p-2)\Delta_{\infty}^{N}v) & \text{ in }Q,\\
 v(0,0)=1,&
\end{cases}
\end{split}
\]
where $Q:=B_4(0)\times (-4C, 4C)$. Observe that  $Q$ is obtained of $B_{4r}(x_0)\times (t_0-4\theta, t_0+4\theta)$ in this rescaling.
Now it suffices to show that there are $\theta_0,\mu_0>0$ so that $\inf_{x\in B_1(0)}v(x,\theta_0)\ge \mu_0$.

\noindent {Step 1: Oscillation estimate.}
To this end, we consider the cylinders
$Q_{\rho,\rho^q}:=Q_{\rho,\rho^q}(0,0):=B_{\rho}(0)\times (-\rho^q,0), \rho\in (0,1),$
and
\[
\begin{split}
M(\rho):=\begin{cases}
\sup_{Q_{\rho,\rho^q}} v, & \rho\in (0,1)\\
1, & \rho=0
\end{cases}, \qquad N(\rho):=(1-\rho)^{-\beta}
\end{split}
\]
where $\beta>1$ will be fixed later. Take $\rho_0\in[0,1)$ be the largest root for the equation $M(\rho)=N(\rho)$. Such a root exists since
\[
\begin{split}
M(0)=1=N(0),\qquad \lim_{\rho\to 1}M(\rho)<\infty,\ \lim_{\rho\to 1}N(\rho)=\infty,
\end{split}
\]
and the functions are continuous on $[0,1)$.
In particular,
\begin{equation}
\label{eq:sup-smaller}
\begin{split}
\sup_{Q_{\rho,\rho^q}} v\le N(\rho),\text{ for all }1>\rho>\rho_0.
\end{split}
\end{equation}
By the continuity of $v$, there is within $\ol Q_{\rho_0,\rho_0^q}$ a point $(x',t')$ such that
\begin{equation}
\label{eq:value-at-x0}
\begin{split}
v(x',t')=\sup_{Q_{\rho_0,\rho_0^q}} v=N(\rho_0)=(1-\rho_0)^{-\beta}.
\end{split}
\end{equation}
Set $R=\half (1-\rho_0)$ (i.e.\ $R$ depends on the $\sup v$) and $Q_{R,R^q}(x',t')=B_{R}(x')\times (t'-R^q,t')$. It holds that $
Q_{R,R^q}(x',t')\subset Q_{\half(1+\rho_0),(\half(1+\rho_0))^q}(0,0)$
so that
\begin{equation}
\label{eq:om}
\begin{split}
\sup_{Q_{R,R^q}(x',t')}v \le M(\half(1+\rho_0))\le N(\half(1+\rho_0))=2^{\beta}(1-\rho_0)^{-\beta}=:\om_0>1.
\end{split}
\end{equation}
Since $\om_0>1$, it holds that $a_0R^q=R^q/\om_0^{q-2}<R^q$ so that $Q_{R,a_0R^q}\subset Q_{R,R^q}$ and
\[
\begin{split}
\sup_{Q_{R,a_0R^q}(x',t')}v\le \om_0.
\end{split}
\]
Thus the assumption of Lemma \ref{lem:key-osc} is satisfied, and Corollary \ref{cor:key-osc} is at our disposal. It follows that there is $\hat C>1,\a \in (0,1)$ such that
\[
\begin{split}
\osc_{x\in B_{r}(x')} v(x,t')\le \hat C\om_0\Big(\frac{r}{R} \Big)^{\a}.
\end{split}
\]

Let $r=\delta R$, $x\in B_{\delta R}(x')$ and observe by the previous estimate together with (\ref{eq:value-at-x0}) that for small enough $\delta>0$ it holds that
\begin{equation}
\label{eq:small-ball-lower-bound}
\begin{split}
v(x,t')&\ge v(x',t')- \hat C \Big(\frac{\delta R}{R}\Big)^{\alpha}2^{\beta}(1-\rho_0)^{-\beta}\\
&\ge (1-\hat C \delta^{\alpha}2^{\beta})(1-\rho_0)^{-\beta}\ge \half (1-\rho_0)^{-\beta}=:\eta.
\end{split}
\end{equation}
Observe that the choice of $\delta>0$ only depends on $n,p,q,\beta$.

\noindent {Step 2: Expansion of positivity.}
Next we use the radial solution (\ref{eq:barenblatt-q-d}) to expand the positivity by using the comparison principle. Without loss of generality we may assume that $(x',t')=(0,0)$.
 We use the Barenblatt type solution from (\ref{eq:barenblatt-q-d}) i.e.\ $\mathcal B^{q,d}(x,\frac{p-1}{q-1}t)$ where
\[
\begin{split}
 \mathcal B^{q,d}(x,t)&=t^{-d/\lambda}\Bigg(C-\frac{q-2}{q}\lambda^{\frac{1}{1-q}}\Big(\frac{\abs x}{t^{1/\lambda}}\Big)^{\frac{q}{q-1}}\Bigg)_+^{\frac{q-1}{q-2}}
\\
&=t^{-d/\lambda}\Big(\frac{q-2}{q}\lambda^{\frac{1}{1-q}}\Big)^{\frac{q-1}{q-2}}\Bigg(C-\Big(\frac{\abs x}{t^{1/\lambda}}\Big)^{\frac{q}{q-1}}\Bigg)_+^{\frac{q-1}{q-2}}\\
&=b t^{-d/\lambda}\Bigg(C-\Big(\frac{\abs x}{t^{1/\lambda}}\Big)^{\frac{q}{q-1}}\Bigg)_+^{\frac{q-1}{q-2}},
\end{split}
\]
with $b=\Big(\frac{q-2}{q}\lambda^{\frac{1}{1-q}}\Big)^{\frac{q-1}{q-2}}$ and $C$ varies from line to line. The scaling $u(x,t/a^{q-2})/a$ preserves the solution. Thus choosing $a=b\nu^{-1}$ we see that
$$
\mathcal B^{q,d}\bigg(x,\frac{p-1}{q-1}\frac{t}{(b\nu^{-1})^{q-2}}\bigg)\frac1{b\nu^{-1}}=:\mathcal B^{q,d}\bigg(x,S(t)\bigg)\frac1{b\nu^{-1}}=\nu S^{-d/\lambda}(t)\Big(1-\Big(\frac{\abs x}{S^{1/\lambda}(t)}\Big)^{\frac{q}{q-1}}\Big)_+^{\frac{q-1}{q-2}},
$$
where  $\nu$ is to be chosen later and $C$ was chosen in a suitable manner, is a solution. Here
\[
\begin{split}
S(t):=\frac{p-1}{q-1}\frac{t}{(b\nu^{-1})^{q-2}}.
\end{split}
\]
The solution is also preserved by a translation of the $t$ variable so that we can consider the solution
\[
\begin{split}
\tilde {\mathcal B}(x,t):=\nu (S(t)+\tau_0)^{-d/\lambda}\Big(1-\Big(\frac{\abs x}{(S(t)+\tau_0)^{1/\lambda}}\Big)^{\frac{q}{q-1}}\Big)_+^{\frac{q-1}{q-2}}.
\end{split}
\]

We intend to select $\nu$ and $\tau_0$ so that
\begin{equation}
\label{eq:req1}
\begin{split}
\spt \tilde {\mathcal B}(\cdot, 0)&\subset \ol B_{\delta R}(0),\\
\tilde {\mathcal B}(\cdot, 0)\le \eta& \text{ in } \ol B_{\delta R}(0),
\end{split}
\end{equation}
where $\eta$ is as in (\ref{eq:small-ball-lower-bound}). To guarantee the first requirement, it suffices to choose $\tau_0$ so that
\[
\begin{split}
(S(0)+\tau_0)^{1/\lambda}=(0+\tau_0)^{1/\lambda}=\delta R,
\end{split}
\]
i.e.\ we can choose $\tau_0=(\delta R)^{\lambda}$.
Then we select $\nu$ so that the second requirement in (\ref{eq:req1}) is satisfied. To guarantee this, it suffices to choose $\nu$ so that
\[
\begin{split}
\tilde {\mathcal B}(0, 0)=\nu (S(0)+\tau_0)^{-d/\lambda}=\nu ((\delta R)^{\lambda})^{-d/\lambda}\le \eta
\end{split}
\]
 i.e.\ we can choose $\nu=\eta (\delta R)^{d}$.

Next we fix $\beta=d$ for the $\beta$ in (\ref{eq:small-ball-lower-bound}). Then we solve for the largest time $\tilde t$ with $\abs x=2$ for which $\tilde {\mathcal B}(x,\tilde t)=0$  from
\[
\begin{split}
1-\Big(\frac{2}{(S(\tilde t)+(\delta R)^{\lambda})^{1/\lambda}}\Big)^{\frac{q}{q-1}}=0.
\end{split}
\]
This  gives
\[
\begin{split}
2^{\lambda}-(\delta R)^{\lambda}=S(\tilde t)=\frac{p-1}{q-1}\frac{\tilde t}{(b\nu^{-1})^{q-2}}.
\end{split}
\]
In other words,
\[
\begin{split}
\tilde t&=(2^{\lambda}-(\delta R)^{\lambda})\frac{q-1}{p-1}(b\nu^{-1})^{q-2}\\
&=(2^{\lambda}-(\delta R)^{\lambda})\frac{q-1}{p-1}\big({(q-2)\lambda^{\frac{1}{1-q}}} q^{-1}\big)^{q-1}\big(\eta (\delta R)^{d}\big)^{2-q}\\
&\ge (2^{\lambda}-1)\frac{q-1}{p-1}\big({(q-2)\lambda^{\frac{1}{1-q}}}q^{-1}  \big)^{q-1}\big(\half (1-\rho_0)^{-d} (\delta \half (1-\rho_0))^{d}\big)^{2-q}\\
&=(2^{\lambda}-1)\frac{q-1}{p-1}\big({(q-2)\lambda^{\frac{1}{1-q}}}q^{-1}  \big)^{q-1}\big(\delta^{d} (\half)^{d+1} \big)^{2-q}.
\end{split}
\]
Above we recalled that $R=\half (1-\rho_0),$  $\delta R<\half$ and $\eta=\half(1-\rho_0)^{-d}$, to see that there is a uniform lower bound for $\tilde t$.
With the choices of the parameters made above, we have
\[
\begin{split}
u\ge \tilde {\mathcal B}\text{ on } \partial_p (B_2(0)\times (0,\tilde t)).
\end{split}
\]
Setting $\theta_0:=\tilde t$, the comparison principle then implies that there is a uniform lower bound $\mu_0$ such that  $\inf_{x\in B_{1}(0)}u(x,\theta_0)\ge \mu_0>0$, so that we have found  $\theta_0,\mu_0$ as intended at the beginning of the proof.

The case $q\le 2$ is rather similar, see \cite{dibenedetto93}.
\qed

The above Harnack's inequality implies the following corollary where $\theta$  is prescribed independently of the solution. The proofs are similar to those in \cite{dibenedetto93}, Theorem 2.2 and Corollary 2.1 on p.158--159.
\begin{corollary}
Let $u\ge 0$ be a viscosity solution to (\ref{eq:qp}) in $Q_{1,1}$ and $q>2$. Then there exists  $C=C(n,p,q)$ such that for all $(x_0,t_0)\in Q_{1,1}$ and for all $r, \theta$ such that $B_{4r}(x_0)\times (t_0-4\theta, t_0+4\theta)\subset Q_{1,1}$ it holds that
\[
\begin{split}
u(x_0,t_0)\le C\Bigg\{ \Big(\frac{r^q}{\theta}\Big)^{\frac1{q-2}}+\Big(\frac{\theta}{r^q}\Big)^{\frac{d}{q}}\Big[ \inf_{y\in B_r(x_0)}u(y,t_0+\theta) \Big]^{\frac{\lambda}{q}} \Bigg\},
\end{split}
\]
where $\lambda=d(q-2)+q$. Under the same conditions, it also holds that
\[
\begin{split}
\kint_{B_r(x_0)} u(x,t_0)\ud x\le C\Bigg\{ \Big(\frac{r^q}{\theta}\Big)^{\frac1{q-2}}+\Big(\frac{\theta}{r^q}\Big)^{\frac{d}{q}}\Big[ u(x_0,t_0+\theta) \Big]^{\frac{\lambda}{q}} \Bigg\}.
\end{split}
\]
\end{corollary}

\section{Below the range condition}

The theory we have displayed solves a number of basic questions for the general $(q,p)$ equation \eqref{eq:qp}, and at the same time leads to some open questions.

Next we discuss the role of the range condition
\begin{equation}
\label{eq:pre-range2}
 2n < q(n-1)+2p,
\end{equation}
in $\Rn$.
 What happens for exponents below this range?
Let us take $q=p$ to simplify matters. The range condition reads then
$$
p>p_c=\frac{2n}{n+1}\,.
$$
It is well-known that for $1<p<p_c$, the $p$-Laplacian theory undergoes a large number of differences with respect to the  case $p>p_c$. One of them is the existence of solutions that extinguish identically in finite time, cf.\ \cite{vazquez06}.

The property of extinction in finite time has been studied in great detail for the Porous Medium Equation, $\partial_t u=\Delta u^m$, PME-$m$, and many results are described in the last reference. In particular the critical exponent is $m_c=(n-2)/n$ for $n\ge 3$. Many types of solutions with finite time extinction can be constructed for $0<m<m_c$, and a number of them are reported in  \cite{vazquez06}.

On the other hand, there is a transformation that maps radial solutions of the PME-$m$ in space dimension $n$ into radial solutions of the $p$-Laplacian equation in a different dimension,
$$
n_1=(n-2)\frac{m+1}{2m},
$$
provided that $p=m+1$. Note that both dimensions need not be integers,  all calculations are made for weighted 1-D equations. This surprising result has been established in \cite{ISV08} and the solution of the $p$-Laplacian equation that is produced is a function $u_1(r',t)$ given by
$$
\partial_{r'} u_1(r',t)=C r^{2/(m+1)}u(r,t), \quad r'=r^{2m/(m+1)}\,,
$$
where $u(r,t) $ is a solution of the PME-$m$,  and $C$ is an inessential constant. Using the PME critical value $m_c=(n-2)/n$ and working out the details of the transformation, we get the corresponding critical value for the $p$-Laplacian equation  $p_c=2n_1/(n_1+1)$. In this way lots of extinguishing solutions can be obtained for the $p$-Laplacian equation if $1<p<p_c$.

 For radial solutions the general $(q,p)$ equation \eqref{eq:qp} reduces to the standard $q$-Laplacian in the fictitious dimension $d$, so we conclude that our range condition \eqref{eq:pre-range2}
marks indeed the border with the possible occurrence of extinction.

In the case of non-radial solutions, all these equations are not equivalent and the theory has to be carefully developed.


\appendix

\section{Equivalence theorem}\label{sec.equiv}

Here we prove the equivalence of viscosity and weak solutions stated in Theorem \ref{thm:equiv}. The proof is divided into two propositions, Propositions \ref{prop:weak-is-visc} and \ref{prop:visc-is-weak}.

First, we recall uniqueness and comparison results for weak solutions.
\begin{lemma}
\label{lem:uniqueness}
Let $u$ and $v$ be two weak solutions according to Definition \ref{def:1D-weak} with $u,v,u_r,v_r\in C([-R,R]\times [0,T))$, $R<\infty$. If $u=v$ on $\partial_p((-R,R)\times(0,T))$, then $u=v$ in $(-R,R)\times(0,T)$.
\end{lemma}
\noindent {\sc Proof.}
Let $u$ and $v$ be two weak solutions.
We test the weak formulations to $u$ and $v$ with
\[
\begin{split}
\phi(r,t)=\chi_h^{0,t_1}(t)(u(r,t)-v(r,t))
\end{split}
\]
with
\[
\begin{split}
\chi_h(t):=\chi_h^{0,t_1}(t)=
\begin{cases}
0 & t\le h ,\\
(t-h)/h, & h<t\le 2h, \\
1,& 2h<t\le t_1-2h,\\
(-t+t_1-h)/h,& t_1-2h<t\le t_1-h,\\
0,& t_1-h<t.
\end{cases}
\end{split}
\]
By a standard approximation argument that we omit, this is admissible. We subtract the weak formulations to obtain
\begin{equation}
\label{eq:subtr-eq}
\begin{split}
\frac{p-1}{q-1}\int_{(-R,R)\times(0,T)} (\abs{v_{r}}^{q-2}v_{r}-\abs{v_{r}}^{q-2}v_{r})\cdot \phi_{r} \ud z=\int_{(-R,R)\times(0,T)} (u-v)\parts{\phi}{t} \ud z.
\end{split}
\end{equation}

 We estimate
\[
\begin{split}
\int_{(-R,R)\times(0,T)} &(u-v) \parts{\phi}{t} \ud z\\
&=\int_{(-R,R)\times(0,T)} (u-v) \parts{(\chi_h (u-v))}{t} \ud z\\
&=\int_{(-R,R)\times(0,T)} (u-v) \Big(\parts{\chi_h}{t}(u-v)+ \chi_h \parts{(u-v)}{t}\Big) \ud z\\
&=\int_{(-R,R)\times(0,T)} (u-v)^2 \parts{\chi_h}{t}\ud z+\int_{(-R,R)\times(0,T)}  \chi_h \half \parts{(u-v)^2}{t} \ud z.
\end{split}
\]
Then we integrate by parts and pass to the limit
\[
\begin{split}
\int_{(-R,R)\times(0,T)} &(u-v)^2 \parts{\chi_h}{t}\ud z-\half \int_{(-R,R)\times(0,T)}  \parts{\chi_h}{t} (u-v)^2 \ud z\\
&=\half \int_{(-R,R)\times(0,T)}  \parts{\chi_h}{t} (u-v)^2 \ud z\\
&=\frac1{2h} \int_h^{2h} \int_{(-R,R)}   (u-v)^2\ud z-\frac1{2h} \int_{t_1-2h}^{t_1-h} \int_{(-R,R)}   (u-v)^2\ud z\\
&\stackrel{\text{$h\to 0$}}{=} 0-\half \int_{(-R,R)}   (u(r,t_1)-v(r,t_1))^2\abs{r}^{d-1}\ud r,
\end{split}
\]
where at the last step we used the initial condition.

By using a well-known algebraic inequality on the right hand side of (\ref{eq:subtr-eq}), and combining the estimates, we obtain with $C>0$
\[
\begin{split}
0&\ge \half \int_{(-R,R)}   (u(r,t_1)-v(r,t_1))^2 \abs{r}^{d-1}\ud r+C\int_{(-R,R)\times(0,T)}\abs{v_{r}-v_{r}}^q \ud z.
\end{split}
\]
Thus since the weight $\abs{r}^{d-1}>0$ whenever $r\neq 0$, we get $u=v$ in $(-R,R)\times(0,T)$.
\qed

\

The above proof also immediately gives $L^2$-contraction property. Indeed, if the initial values are $u_0,v_0$ and the lateral boundary values are the same, the last inequality in the above proof reads as
\[
\begin{split}
 \half &\int_{(-R,R)}   (u_0(r)-v_0(r))^2 \abs{r}^{d-1}\ud r\\
 &\ge \half \int_{(-R,R)}   (u(r,t_1)-v(r,t_1))^2 \abs{r}^{d-1}\ud r+C\int_{(-R,R)\times(0,T)}\abs{v_{r}-v_{r}}^q \ud z\\
 &\ge  \half \int_{(-R,R)}   (u(r,t_1)-v(r,t_1))^2 \abs{r}^{d-1}\ud r.
\end{split}
\]
Moreover, if we test with $\phi(r,t)=\chi_h^{0,t_1}(t)\abs{u(r,t)-v(r,t)}^{m-1}(u(r,t)-v(r,t)),\ 1<m<\infty$ instead, then a similar computation as above gives
\[
\begin{split}
 \frac{1}{m} &\int_{(-R,R)}   (u_0(r)-v_0(r))^m \abs{r}^{d-1}\ud r\\
 &\ge \frac{1}{m} \int_{(-R,R)}   (u(r,t_1)-v(r,t_1))^m \abs{r}^{d-1}\ud r+C\int_{(-R,R)\times(0,T)}\abs{u-v}^{m-2}\abs{v_{r}-v_{r}}^q \ud z\\
 &\ge  \frac{1}{m} \int_{(-R,R)}   (u(r,t_1)-v(r,t_1))^m \abs{r}^{d-1}\ud r.
\end{split}
\]
To make this rigorous one would have to mollify in time.

A similar proof to the uniqueness also gives a comparison principle.
\begin{lemma}
\label{lem:comparison}
Let $u$ be a weak subsolution and $v$ a weak supersolution  according to Definition \ref{def:1D-weak}  with $u,v,u_r,v_r\in C([-R,R]\times [0,T))$, $R<\infty$. If $u\le v$ on $\partial_p((-R,R)\times(0,T))$, then $u\le v$ in $(-R,R)\times(0,T)$.
\end{lemma}

\begin{proposition}
\label{prop:weak-is-visc}
Let $u\in C(Q_T)$, $Q_T=B_R\times(0,T), B_R\subset \R^n,  0<R\le \infty,$ be a continuous radial function, and $q>1,p>1$.
If  $v(r,t):=u(re_1,t),\ r\in (-R,R),$ is 1-dimensional  weak solution to (\ref{eq:1D-div-form})
according to Definition \ref{def:1D-weak}, then $u$ is a viscosity solution to (\ref{eq:qp}) in $n$-dimensions.
\end{proposition}
\noindent {\sc Proof.}
Since the case of sub- and supersolutions is analogous, thriving for a contradiction, we may assume that there is an admissible (according to Definition \ref{def:admissible})
test function $\vp\in C^2$  touching $u$ from below at $(x_0,t_0)\in Q_T$ and one of the two cases holds
\[
\begin{split}
\begin{cases}
\vp_t(x_0,t_0)-F(D\vp(x_0,t_0),D^2\vp(x_0,t_0))<0,& \text{ if }D\vp(x_0,t_0)\neq	 0 \\
 \vp_t(x_0,t_0)<0,& \text{ if }D\vp(x_0,t_0)=0.
\end{cases}
\end{split}
\]

Consider first the case $D\vp(x_0,t_0)\neq 0$, $x_0\neq 0$. With the usual abuse of notation, we keep using $\vp$ also when in spherical coordinates, and $r>0$. Then it holds recalling $u$ is radial that
\begin{equation}
\label{eq:CP-from-below}
\begin{split}
\vp_t<F(D\vp,D^2\vp)&=\abs{\vp_{r}}^{q-2}\Big(\vp_{rr}+\frac{n-1}{r}\vp_r+\frac1{r^2}\Delta_{S^{n-1}}\vp+(p-2)\vp_{rr}\Big)\\
&\le \abs{\vp_{r}}^{q-2}\Big((p-1)\vp_{rr}+\frac{n-1}{r}\vp_r\Big)\\
&=\frac{p-1}{q-1}\abs{\vp_{r}}^{q-2}\Big(\frac{d-1}{r}\vp_r+(q-1)\vp_{rr}\Big)
\end{split}
\end{equation}
where $\Delta_{S^{n-1}}$ is the Laplace-Beltrami operator on $n-1$-sphere. Next set $\phi(r,t):=\vp(r\frac{{x_0}}{\abs {x_0}},t)$.  Since $\phi\in C^2$, it is a 1-dimensional weak subsolution according to Definition \ref{def:1D-weak} in some cylinder $Q^1_{\delta,\delta}:=Q^1_{\delta,\delta}(x_0,t_0):=(\abs{{x_0}}-\delta,\abs{{x_0}}+\delta)\times (-\delta+t_0,t_0)$  by the computation (\ref{eq:1-divergence-form}) and (\ref{eq:CP-from-below}). Then contradiction follows by a standard argument, i.e.\ adding a constant $m>0$ small enough such that $\emptyset\neq \{  \phi+m>u \}\Subset Q^1_{\delta,\delta}$ and $\phi_r\neq 0$ in $\{  \phi+m>u \}$. Since $ \phi+m\in C^2$ is also a weak subsolution and $u$ is a weak solution, we arrive at the contradiction recalling the comparison principle, Lemma \ref{lem:comparison}.

\sloppy Consider then the case $D\vp(x_0,t_0)=0,\ x_0\neq 0$, and let us assume, in the search of a contradiction, that $u_t(x_0,t_0)=\vp_t(x_0,t_0)<0$. Moreover, by \cite[Remark 2.2.7]{giga06} we may assume that $D\vp(x,t)\neq 0$ whenever $x\neq x_0$, and
\[
\begin{split}
\lim_{x_0\neq x\to x_0,}F(D\vp(x,t_0),D^2\vp(x,t_0))=0.
\end{split}
 \]
By this and the counter assumption, denoting $Q_{\delta,\delta}=B_{\delta}(x_0)\times (t_0-\delta, t_0)$, we have
\[
\begin{split}
\vp_t<F(D\vp,D^2\vp)
\end{split}
\]
in $\{ (x,t)\in Q_{\delta,\delta} \,:\,x\neq x_0\}$ for small enough $\delta>0$. Since $u$ is radial and thus $\Delta_{S^{n-1}}\vp/r^2\le 0$, it follows by continuity of $\Delta_{S^{n-1}}\vp$ and  by combining the calculations (\ref{eq:1-divergence-form}) and (\ref{eq:CP-from-below}) that
\[
\begin{split}
\phi_t-\frac{p-1}{q-1}\Big(\abs{\phi_r}^{q-2}\phi_r \abs  r^{d-1}\Big)_r \abs r^{1-d}<0
\end{split}
\]
in $\{ (r,s)\in Q^1_{\delta,\delta}\,:\,r\neq \abs {x_0}\}$. Without loss of generality, we may take $\delta>0$ small enough so that for the notational convenience $r>0$. Using this with $\eta\in C^\infty_0(Q^1_{\delta,\delta}),\,\eta\ge 0$, we obtain
\[
\begin{split}
&\int_{Q^1_{\delta,\delta}} \abs{\phi_r}^{q-2}\phi_r r^{d-1}  \eta_r \ud r\ud t=\lim_{\rho\to 0} \int_{Q^1_{\delta,\delta}\setminus \{(r,t)\,:\,\abs{\abs {x_0}-r}\le \rho \}} \abs{\phi_r}^{q-2}\phi_r r^{d-1}  \eta_r \ud r\ud t\\
&=\lim_{\rho\to 0} \Big\{-\int_{Q^1_{\delta,\delta}\setminus \{(r,t)\,:\,\abs{\abs {x_0}-r}\le \rho \}} \Big(\abs{\phi_r}^{q-2}\phi_r r^{d-1}\Big)_r \eta \ud r\ud t-\int_{t_0-\delta}^{t_0}\Big[\abs{\phi_r}^{q-2}\phi_r r^{d-1}\eta\Big]^{\abs {x_0}+\rho}_{\abs {x_0}-\rho}\ud t\Big\}\\
&\le  -\frac{q-1}{p-1}\int_{Q^1_{\delta,\delta}}\phi_t r^{d-1} \eta \ud r\ud t= \frac{q-1}{p-1}\int_{Q^1_{\delta,\delta}}\phi r^{d-1} \eta_t \ud r\ud t,
\end{split}
\]
where in the first step we used the dominated convergence theorem and the fact that $q-1>0$. This again implies that $\phi$ is a weak subsolution, and the contradiction is obtained similarly as in the first case.

Finally, consider  the case $D\vp(x_0,t_0)=0, x_0=0$ (the weak solution has $v_r(0,t)=0$ so $D\vp(x_0,t_0)\neq 0$ does not occur), and observe that the argument in the previous case only utilized the equation at $x\neq x_0$. Moreover, our test function in this case can be taken to be of the form $\vp(x,t)=f(\abs{x})+g(t)$, \cite[Remark 2.2.7]{giga06}, which is a radial $C^{2}$-function in $\Rn$.  Thus it holds that $\Delta_{S^{n-1}} \vp=0$ outside the origin, and thus the computation similar to (\ref{eq:CP-from-below}), recalling $d-1>0$, still holds. The contradiction then follows similarly as before.
\qed

Next we show that a radial viscosity solution is a 1-dimensional weak solution.

\begin{proposition}
\label{prop:visc-is-weak}
Let $u\in C(Q_T)$, $Q_T=B_R\times(0,T), B_R\subset \R^n,   0<R\le \infty$ be a  continuous radial function, and $q>1$.
 Then if $u$ is a viscosity solution to (\ref{eq:qp}) in $n$-dimensions, it follows that  $v(r,t):=u(re_1,t),\ r\in (-R,R),$ is 1-dimensional  weak solution to (\ref{eq:1D-div-form})
according to Definition \ref{def:1D-weak}.
\end{proposition}

\noindent {\sc Proof.}
This is a parabolic version of the proof in \cite{julinj12}. Without a loss of generality, we may assume in the proof that $R<\infty$ and work out the proof in $B_R\times (0,T)$ even  in the case $\Rn\times (0,T)$. If $u$ is a weak solution in $B_R\times (0,T)$ for all $R<\infty$, then it is a weak solution in  $\Rn\times (0,T)$.

First, suppose that $q> 2$.
We will prove the weak supersolution property; the proof of the subsolution property is similar. To be more precise, we  show that $v(r,t):=u(re_1,t)$ satisfies
\begin{equation}
\label{eq:aim-super}
\int_{Q_T^1} v \phi_t \ud z\le \frac{p-1}{q-1}\int_{Q_T^1} \abs{v_{r}}^{q-2}v_{r} \phi_{r} \ud z
\end{equation}
where $\phi\in C^{\infty}_0(Q_T^1), \phi\ge 0,$ where $Q_T^1:=(-R,R)\times(0,T)$.  The $C^1$-conditions in the definition are immediately satisfied:  $v(\cdot,t)$ is  $C^1$-function  and $v_r(0,t)=0$ because $u$ is $C^{1}$ by \cite{imbertjs}, and $u$ is radial.

\medskip

\noindent{Step 1: Regularization.} Let us continue by showing that the inf-convolution $u_\eps$ of $u$,
\begin{equation}
\label{eq:inf-conv}
u_{\eps}(x,t):=\underset{(y,s)\in Q_T}{\inf}\left( u(y,s)+\dfrac{|x-y|^{2}+\abs{t-s}^2}{2\eps}\right),
\end{equation}
is a weak supersolution in
$$Q_{\eps}=\left\{(x,t)\ :\  \dist((x,t), \partial Q_T)> (2\eps \osc_{Q_T}{u})^{1/2}\right\}.$$

First, it holds that $u_{\eps}$ is a semiconcave viscosity supersolution to (\ref{eq:qp}). The Sobolev derivatives $\partial_t u_{\eps}, D u_{\eps}$ exist and belong to $L^{\infty}_{\text{loc}}(Q_{\eps})$. Moreover, $u_\eps$  is semiconcave and twice differentiable a.e.\ and satisfies
\begin{align} \nonumber
 \partial_t u_{\eps}\ge |Du_{\eps}|^{q-2}\left(\Delta u_{\eps}+(p-2) D^2u_{\eps}\dfrac{Du_{\eps}}{|Du_{\eps}|}\cdot \dfrac{Du_{\eps}}{|D u_{\eps}|}\right)
\end{align}
a.e.\ in $Q_{\eps}$, and $u_{\eps}$ is still radial. For the properties of parabolic infimal convolutions, see for example \cite{lindqvist12}. Also observe that since $q>2$, the interpretation of the right hand side is  clear also if $Du_{\eps}(x,t)=0$.  It follows that in radial coordinates it holds similarly as before
\begin{equation}
\begin{split}
\label{eq:eq-for-inf}
 \partial_t u_\eps&\ge F(Du_{\eps},D^2u_\eps)\\
&=\abs{(u_\eps)_{r}}^{q-2}\Big((u_\eps)_{rr}+\frac{n-1}{r}(u_\eps)_r+\frac1{r^2}\Delta_{S^{n-1}}(u_\eps)+(p-2)(u_\eps)_{rr}\Big)\\
&= \abs{(u_\eps)_{r}}^{q-2}\Big((p-1)(u_\eps)_{rr}+\frac{n-1}{r}(u_\eps)_r\Big)\\
&=\frac{p-1}{q-1}\abs{(u_\eps)_{r}}^{q-2}\Big((q-1)(u_\eps)_{rr}+\frac{d-1}{r}(u_\eps)_r\Big)\\
&=\frac{p-1}{q-1}(\abs{(u_\eps)_{r}}^{q-2}(u_\eps)_{r}\abs{r}^{d-1})_r \abs{r}^{1-d},
\end{split}
\end{equation}
 a.e. in $Q_{\eps}$.  In particular, we may assume that $r\neq 0$, since $\{ (r,t) \,:\, r=0 \}$ is of measure zero.

Since $u_{\eps}$ is semiconcave i.e.\ the function $(x,t)\mapsto u_{\eps}(x,t)-\dfrac{1}{2\eps}(|x|^2+t^2)$ is concave  in $Q_{\eps}$, we can approximate it by a sequence $(\vp_j)$ of smooth concave radial functions by using the standard mollification. Denoting $v_{\eps, j}(r,t):=\vp_j(re_1,t)+\dfrac{1}{2\eps}(|r|^2+t^2)$, we can integrate by parts to obtain
\begin{equation}\label{hep1}
-\int_{Q_T^1}(\abs{(v_{\eps, j})_{r}}^{q-2}(v_{\eps, j})_{r}\abs{r}^{d-1})_r \phi\, dr\ud t=\int_{Q_T^1} \abs{(v_{\eps, j})_{r}}^{q-2}(v_{\eps, j})_{r}\abs{r}^{d-1} \phi_r dr\ud t,
\end{equation}
\sloppy
 for any nonnegative $\phi\in C_0^{\infty}(Q_{\eps}^1)$ where $Q_{\eps}^1=\{(r,t)\in Q_T^1\,:\,\dist((r,t), \partial Q_T^1)> (2\eps \osc_{Q_T^1}{u})^{1/2}\}$.
Since $Du_{\eps}$ and thus  $(v_{\eps})_r$  are in $L^{\infty}_{\text{loc}}$, the dominated convergence theorem implies
\begin{equation}\label{hep2}
\underset{j\to\infty}{\lim} \int_{Q_T^1} \abs{(v_{\eps, j})_{r}}^{q-2}(v_{\eps, j})_{r}\abs{r}^{d-1} \phi_r dr\ud t
=\int_{Q_T^1} \abs{(v_{\eps})_{r}}^{q-2}(v_{\eps})_{r}\abs{r}^{d-1} \phi_r dr\ud t.
\end{equation}

Next, by concavity of $\vp_j$ we have $(v_{\eps, j})_{rr}\leq \frac1{\eps}$ and thus  by the local boundedness of $(v_{\eps,j})_r$, we get
\[
\begin{split}
-(\abs{(v_{\eps, j})_{r}}^{q-2}(v_{\eps, j})_{r}\abs{r}^{d-1})_r&=-\abs{(v_{\eps, j})_{r}}^{q-2}\Big(\frac{d-1}{r}(v_{\eps, j})_r+(q-1)(v_{\eps, j})_{rr}\Big)\abs{r}^{d-1}\\
&\ge -C^{q-2}( C(d-1) \abs{r}^{d-2}+(q-1)\abs{r}^{d-1}/\eps ).
\end{split}
\]
Since $d>1$, this is an integrable lower bound needed for Fatou's theorem.
Applying Fatou's theorem, we obtain
\begin{equation}
\label{hep3}
\begin{split}
\underset{j\to\infty}{\liminf}&\int_{Q_T^1} -(\abs{(v_{\eps, j})_{r}}^{q-2}(v_{\eps, j})_{r}\abs{r}^{d-1})_r\phi \ud r\ud t \\
&\geq \int_{Q_T^1} \underset{j\to\infty}{\liminf} -(\abs{(v_{\eps, j})_{r}}^{q-2}(v_{\eps, j})_{r}\abs{r}^{d-1})_r \phi \ud r \ud t.
\end{split}
\end{equation}
Since
$$\underset{j\to\infty}{\liminf} -(\abs{(v_{\eps, j})_{r}}^{q-2}(v_{\eps, j})_{r}\abs{r}^{d-1})_r= -(\abs{(v_{\eps})_{r}}^{q-2}(v_{\eps})_{r}\abs{r}^{d-1})_r$$
almost everywhere, by using \eqref{hep1}, \eqref{hep2} and \eqref{hep3} we obtain
\begin{align}
\label{eq:super-weak}
\frac{p-1}{q-1}\int_{Q_T^1} \abs{(v_{\eps})_{r}}^{q-2}(v_{\eps})_{r}&\abs{r}^{d-1} \phi_r \ud r\ud t \geq \frac{p-1}{q-1} \int_{Q_T^1}-(\abs{(v_{\eps})_{r}}^{q-2}(v_{\eps})_{r}\abs{r}^{d-1})_r\phi\, \ud r\ud t\nonumber\\
& \geq -\int_{Q_T^1} \abs{r}^{d-1}\partial_t v_{\eps}\phi \ud r\ud t=\int_{Q_T^1} \abs{r}^{d-1} v_{\eps}\partial_t \phi \ud r\ud t,
\end{align}
where the last inequality follows from (\ref{eq:eq-for-inf}).
Thus we have accomplished (\ref{eq:aim-super}) but so far only for $v_{\eps}$.

\medskip

\noindent{Step 2: Passing to the limit in the regularization.}
First choose  cylindrical domains $Q''\Subset  Q'\Subset Q_{\eps}^1$. We start by showing that $(v_\eps)_r$ is uniformly bounded in the weighted $L^q(Q'')$. Take a cut-off function  $\xi:Q_{\eps}^1\rightarrow [0,1],\ \xi\in C^{\infty}_0(Q')$ such that $\xi\equiv 1$ on $Q''$. Choose the test function $\phi=(M-v_{\eps})\xi^q$ in (\ref{eq:super-weak}), where $M=\osc_{Q'}|u_\eps|$. Since the test function is Lipschitz by the properties of infimal convolution and compactly supported, this is an admissible test function after an approximation argument. We have
\begin{equation}
\label{eq:time-part}
\begin{split}
\int_{Q_T^1} &\abs{r}^{d-1} v_{\eps}\partial_t \phi dx\ud t\\
&=\int_{Q_T^1}  \abs{r}^{d-1} v_{\eps}\partial_t ((M-v_{\eps})\xi^q)\ud r\ud t\\
&=\int_{Q_T^1} \abs{r}^{d-1} v_{\eps} (-\partial_tv_{\eps}\xi^q+(M-v_{\eps})\partial_t \xi^q)\ud r\ud t\\
&=\int_{Q_T^1} \abs{r}^{d-1} (-\half \partial_t v_{\eps}^2 \xi^q+v_{\eps}(M-v_{\eps})\partial_t \xi^q)\ud r\ud t\\
&=\int_{Q_T^1} \abs{r}^{d-1} (\half  v_{\eps}^2 \partial_t\xi^q+v_{\eps}(M-v_{\eps}) \partial_t \xi^q)\ud r\ud t\le C(p,n,q, \norm{v_{\eps}}_{L^{\infty}(Q')}).
\end{split}
\end{equation}
For
\[
\begin{split}
\int_{Q_T^1} &\abs{r}^{d-1} \abs{(v_{\eps})_{r}}^{q-2}(v_{\eps})_{r} ((M-v_{\eps})\xi^q)_r \ud r\ud t\\
&=\int_{Q_T^1} -\abs{r}^{d-1}\abs{(v_{\eps})_{r}}^{q}  \xi^q \ud r\ud t+\int_{Q_T^1} \abs{r}^{d-1}\abs{(v_{\eps})_{r}}^{q-2}(v_{\eps})_{r} (M-v_{\eps})(\xi^q)_r \ud r\ud t
\end{split}
\]
we use H\"older's inequality. It follows recalling (\ref{eq:super-weak})--(\ref{eq:time-part}) that
\[
\begin{split}
&\int_{Q_T^1}  \abs{r}^{d-1} |(v_{\eps})_r|^q\ \xi^q \ud r \ud t \\
&\leq q\int_{Q_T^1} \abs{r}^{d-1}\xi^{q-1}|(v_{\eps})_r|^{q-2}(v_{\eps})_r \xi_r  (M-v_{\eps})\, dr \ud t+C(p,n,q, \norm{v_{\eps}}_{L^{\infty}(Q')})\\
&\leq\half\int_{Q_T^1} \abs{r}^{d-1}\xi^q|(v_{\eps})_r|^q\,dr\ud t
\\
&\hspace{3 em} +\ C\int_{Q_T^1} \abs{r}^{d-1} M^q|\xi_r|^{q} \ud r \ud t+C(p,n,q, \norm{v_{\eps}}_{L^{\infty}(Q')}).
\end{split}
\]
Absorbing the first term on the right into the left, it follows that
\begin{align}\label{vaartyu}
\int_{Q_T^1} \abs{r}^{d-1} \xi^q|(v_{\eps})_r|^q\,\ud r \ud t
 &\leq C= C\left(p,n,q,\norm{u}_{L^\infty(Q')}\right).
\end{align}
Hence, $(v_{\eps})_r$ is uniformly bounded with respect to $\eps$ in $L^q(\abs{r}^{d-1}\ud r \ud t, Q')$. It follows that there exists a subsequence such that $(v_{\eps})_r\to \tilde v_r$ weakly in $L^{q}(\abs{r}^{d-1}\ud r \ud t,Q')$. Moreover,  choose a smooth test function such that $\spt \phi\subset Q'$ and observe that by the dominated convergence theorem since $d>1$
\[
\begin{split}
\Big|\int_0^T\int_{-R}^{R} \phi \abs{r}^{d-1}&((v_{\eps})_r-v_r)\ud r \ud t\Big|=\lim_{\delta\to 0}\abs{\int_0^T\int_{(-R,-\delta)\cup(\delta,R)} \phi \abs{r}^{d-1}((v_{\eps})_r-v_r) \ud r \ud t}\\
&=\lim_{\delta\to 0}\bigg|-\int_0^T\int_{(-R,-\delta)\cup(\delta,R)}(\phi_r \abs{r}^{d-1}+\phi(d-1)\abs{r}^{d-2})(v_{\eps}-v) \ud r \ud t\\
&\hspace{1 em}+\int_0^T\Big\{\delta^{d-1} (\phi(v_{\eps}-v))(-\delta,t)-\delta^{d-1} (\phi(v_{\eps}-v))(\delta,t)\Big\}\ud t \bigg|\\
&= \abs{\int_0^T\int_{(-R,R)}(\phi_r \abs{r}^{d-1}+\phi(d-1)\abs{r}^{d-2})(v_{\eps}-v) \ud r \ud t}\\
&\le \norm{v_{\eps}-v}_{L^{\infty}(Q')} \abs{\int_0^T\int_{(-R,R)}(\phi_r \abs{r}^{d-1}+\phi(d-1)\abs{r}^{d-2}) \ud r \ud t}.
\end{split}
\]
Moreover the right hand side converges to zero as $\eps\to 0$, so that $\tilde v_r=v_r$ a.e.\ in $Q'$. Then similarly as in Theorem 5.3 in \cite{kortekp10}, see also \cite{lindqvistm07}, it holds that the pointwise limit $v$ of bounded weak supersolutions $v_{\eps}$ is a weak supersolution.

Then consider the case $1<q\le 2$.

\noindent{Step 1: Regularization.}  Again $u_{\eps}$ denotes the $\inf$-convolution of $u$ but now
\begin{equation}
\label{eq:inf-convolution-q}
\begin{split}
u_{\eps}(x,t):=\underset{(y,s)\in Q_T}{\inf}\left( u(y,s)+\frac{|x-y|^{\hat q}}{\hat q\eps^{\hat q-1}}+\frac{\abs{t-s}^2}{2\eps}\right),
\end{split}
\end{equation}
for $\hat q> q/(q-1)$.
 Then similarly as  before in a.e.\ in $Q_{\eps}\setminus \{ Du_{\eps}=0\}$, where the definition of $Q_\eps$ is modified accordingly (see below), the function $u_{\eps}$ is a viscosity supersolution to $ \partial_t u_{\eps}\ge F(Du_{\eps},D^2u_\eps)$, and  we have
\[
\begin{split}
\partial_t u_{\eps}\ge |Du_{\eps}|^{q-2}\left(\Delta u_{\eps}+(p-2) D^2u_{\eps}\dfrac{Du_{\eps}}{|Du_{\eps}|}\cdot \dfrac{Du_{\eps}}{|D u_{\eps}|}\right)= \frac{p-1}{q-1}(\abs{(u_\eps)_{r}}^{q-2}(u_\eps)_{r}\abs{r}^{d-1})_r \abs{r}^{1-d}.
\end{split}
\]
 However, if $Du_{\eps}(x,t)=0$, the meaning of the right hand side is no longer clear. Therefore, we look at the regularized operator in order to integrate by parts using semiconcavity of $u_{\eps}$ and Fatou's theorem, and the dominated convergence theorem with respect to $j$ as in the previous case. We obtain
\begin{equation}
\label{eq:j-fatou}
\begin{split}
\int_{Q_T^1}& ((\abs{(v_{\eps})_{r}}^2+\delta^2)^{\frac{q-2}{2}}(v_{\eps})_{r}\abs{r}^{d-1}) \phi_r  \ud r \ud t\\
&=\lim_{j\to \infty}\int_{Q_T^1} ((\abs{(v_{\eps,j})_{r}}^2+\delta^2)^{\frac{q-2}{2}}(v_{\eps,j})_{r}\abs{r}^{d-1}) \phi_r  \ud r \ud t\\
&\ge-\int_{Q_T^1} \liminf_{j\to \infty}((\abs{(v_{\eps,j})_{r}}^2+\delta^2)^{\frac{q-2}{2}}(v_{\eps,j})_{r}\abs{r}^{d-1})_r \phi  \ud r \ud t\\
&= \int_{Q_T^1} -((\abs{(v_{\eps})_{r}}^2+\delta^2)^{\frac{q-2}{2}}(v_{\eps})_{r}\abs{r}^{d-1})_r \phi  \ud r \ud t,
\end{split}
\end{equation}
for $\phi\in C^\infty_0(Q_{\eps}^1),\ \phi\ge 0$.
 When passing to the limit $\delta\to 0$, we need to justify the convergence also on the right hand side. Let $\hat x$ and $x_{\eps}$ be as in Lemma \ref{lem:inf-conv}. Then by Lemma \ref{lem:inf-conv}
\[
\begin{split}
D u_{\eps}(\hat x, \hat t)&=(\hat x-x_{\eps})\frac{\abs{\hat x-x_{\eps}}^{{\hat q}-2}}{\eps^{{\hat q}-1}},\\
D^{2}u_{\eps}(\hat x,\hat t)&\le \begin{cases}
 (\hat q-1)\frac{\abs{\hat x-x_{\eps}}^{{\hat q}-2}}{\eps^{{\hat q}-1}} I, & \text{if }Du_{\eps}(\hat x, \hat t)\neq 0,\\
0, & \text{if }Du_{\eps}(\hat x, \hat t)= 0.
\end{cases}
\end{split}
\]
Thus denoting $r_{\eps}:=\abs{\hat x-x_{\eps}}$ and $r:=\abs{\hat x}$ we have
\[
\begin{split}
 -&((\abs{(v_{\eps})_{r}}^2+\delta^2)^{\frac{q-2}{2}}(v_{\eps})_{r}\abs{r}^{d-1})_r \\
 &=
 -(\abs{(v_{\eps})_{r}}^2+\delta^2)^{\frac{q-2}{2}}\bigg((q-2)(v_{\eps})_{rr}(v_{\eps})^2_{r}\big((v_{\eps})^2_{r}+\delta^2\big)^{-1}+(v_{\eps})_{rr}+(v_{\eps})_{r}\frac{d-1}{r}\bigg)\abs{r}^{d-1}\\
&\ge -C(q,\hat q,\eps)\bigg( r_{\eps}^{(\hat q-1)(q-1)-1}+ r_{\eps}^{(\hat q-1)(q-1)}\frac{d-1}{r}\bigg)\abs{r}^{d-1}.
\end{split}
\]
Since $\hat q> q/( q-1)$ so that $(\hat q-1)(q-1)-1> 0$, and since $\abs{r}^{d-2}$ is integrable, this gives an integrable lower bound independent of $\delta$.
Thus by Fatou's lemma
\[
\begin{split}
\liminf_{\delta\to 0} &\int_{Q_T^1}  -((\abs{(v_{\eps})_{r}}^2+\delta^2)^{\frac{q-2}{2}}(v_{\eps})_{r}\abs{r}^{d-1})_r \phi \ud r \ud t\\
&\ge\int_{Q_T^1\setminus \{ (v_{\eps})_r =0\}}   \liminf_{\delta\to 0} (-((\abs{(v_{\eps})_{r}}^2+\delta^2)^{\frac{q-2}{2}}(v_{\eps})_{r}\abs{r}^{d-1})_r \phi) \ud r \ud t\\
&= \int_{Q_T^1\setminus \{ (v_{\eps})_r=0\}}   -(\abs{(v_{\eps})_{r}}^{q-2}(v_{\eps})_{r}\abs{r}^{d-1})_r \phi \ud r \ud t\\
&\ge \frac{q-1}{p-1}\int_{Q_T^1\setminus \{ (v_{\eps})_r =0\}}   -\partial_t v_{\eps} \abs{r}^{d-1} \phi \ud r \ud t.
\end{split}
\]
By Lemma \ref{lem:inf-conv}, it follows that  $-\partial_t v_{\eps}\le 0$ in  $\{(v_{\eps})_r=0\}$ so that
\[
\begin{split}
 \int_{Q_T^1\setminus \{ (v_{\eps})_r =0\}}   -\partial_t v_{\eps} \abs{r}^{d-1} \phi \ud r \ud t\ge  \int_{Q_T^1}   -\partial_t v_{\eps} \abs{r}^{d-1} \phi \ud r \ud t= \int_{Q_T^1}    v_{\eps} \abs{r}^{d-1} \partial_t\phi \ud r \ud t.
\end{split}
\]
From this and passing to the limit with $\delta\to 0$ in (\ref{eq:j-fatou}), we obtain
\[
\begin{split}
\frac{p-1}{q-1}\int_{Q_T^1} \abs{(v_{\eps})_{r}}^{q-2}(v_{\eps})_{r}\abs{r}^{d-1} \phi_r \ud r\ud t\ge \int_{Q_T^1}    v_{\eps} \abs{r}^{d-1} \partial_t\phi \ud r \ud t.
\end{split}
\]

\noindent{Step 2: Passing to the limit in the regularization.}  This follows similarly as before. In particular,  techniques similar to those in Theorem 5.3 in \cite{kortekp10} do not utilize a lower bound for $q$ other than $q>1$.
\qed

From Propositions \ref{prop:weak-is-visc} and \ref{prop:visc-is-weak}, Theorem \ref{thm:equiv} immediately follows.

Next we state some auxiliary results used in the proof above. Let $u_{\eps}(x,t)$ be as in (\ref{eq:inf-convolution-q}). Denote
$r(\eps):=(\hat q \eps^{\hat q-1}\osc_{Q_T}u)^{1/\hat q} $, and $t(\eps):=(2 \eps \osc_{Q_T}u)^{1/2}$.
It is well known that for  $(\hat x,\hat t)\in Q_{\eps}=\left\{(x,t)\ :\  B_{r_{\eps}}(\hat x) \Subset B_R, (t- t(\eps),t+t(\eps))\Subset (0,T)\right\}$ there exists $x_{\eps}\in B_{r(\eps)}(\hat x)$ and $t_{\eps}$ with $\abs{\hat t-t_{\eps}}\le t(\eps)$ such that
\[
\begin{split}
u_{\eps}(\hat x,\hat t)=u(x_{\eps},t_{\eps})+\frac{|\hat x-x_{\eps}|^{\hat q}}{\hat q\eps^{\hat q-1}}+\frac{\abs{\hat t-t_{\eps}}^2}{2\eps}.
\end{split}
\]
It holds that $u_{\eps}\nearrow u$
uniformly.
By Alexandrov's theorem, see \cite[Section 6.4]{evansg92}, and semiconcacity of $u_{\eps}$ it follows that $u_{\eps}$ is twice differentiable a.e. We omit the proof of semiconcavity which is well-known, see for example \cite{katzourakis15} or  Lemma A.2 in \cite{julinj12},  and instead derive the following explicit estimates.
\begin{lemma}
\label{lem:inf-conv}
Let $1<q\le 2$, $u_{\eps}$ as in (\ref{eq:inf-convolution-q}), and $x_{\eps}, t_{\eps}$ as above.
Suppose that $u_{\eps}$ is differentiable in time and twice differentiable in space  at $(\hat x,\hat t)$. Then
\begin{enumerate}[(i)]
\item it holds that
\[
\begin{split}
D u_{\eps}(\hat x, \hat t)&=(\hat x-x_{\eps})\frac{\abs{\hat x-x_{\eps}}^{{\hat q}-2}}{\eps^{{\hat q}-1}},\\
D^{2}u_{\eps}(\hat x,\hat t)&\le \begin{cases}
 (\hat q-1)\frac{\abs{\hat x-x_{\eps}}^{{\hat q}-2}}{\eps^{{\hat q}-1}} I & \text{if }Du_{\eps}(\hat x, \hat t)\neq 0, \\
0 & \text{if }Du_{\eps}(\hat x, \hat t)= 0, \text{ and }
\end{cases}
\end{split}
\]
\item if  $Du_{\eps}(\hat x,\hat t)=0$ it follows that $\partial_t u_{\eps}(\hat x,\hat t)\ge 0$.
\end{enumerate}
\end{lemma}
\noindent {\sc Proof.}
 Proof of $(ii)$.
There exists $\vp\in C^{\infty}_0( Q_T)$ touching $u_{\eps}$ at $(\hat x,\hat t)$ from below such that $\partial_t \vp(\hat x,\hat t)=\partial_t u_{\eps}(\hat x, \hat t), D\vp(\hat x,\hat t)=Du_{\eps}(\hat x,\hat t)=0$. Further,
\[
\begin{split}
u(y,s)+\frac{\abs{x-y}^{\hat q}}{{\hat q}\eps^{{\hat q}-1}}+\frac{\abs{t-s}^{2}}{2\eps}-\vp(x,t)&\ge u_{\eps}(x,t)-\vp(x,t)\ge  u_{\eps}(\hat x,\hat t)-\vp(\hat x,\hat t)=0.
\end{split}
\]
Choose $y=x_{\eps},s=t_{\eps}$,  and write
\[
\begin{split}
\vp(x,t)-\Big(\frac{\abs{x-x_{\eps}}^{\hat q}}{{\hat q}\eps^{{\hat q}-1}}+\frac{\abs{t-t_{\eps}}^{2}}{2\eps}\Big)&\le u(x_{\eps},t_{\eps}).
\end{split}
\]
Since
\[
\begin{split}
\vp(\hat x,\hat t)&= u(x_{\eps},t_{\eps})+\frac{\abs{\hat x-x_{\eps}}^{\hat q}}{{\hat q}\eps^{{\hat q}-1}}+\frac{\abs{\hat t-t_{\eps}}^{2}}{2\eps},
\end{split}
\]
it follows that
\[
\begin{split}
(x,t)\mapsto \vp(x,t)-\Big(\frac{\abs{x-x_{\eps}}^{\hat q}}{{\hat q}\eps^{{\hat q}-1}}+\frac{\abs{t-t_{\eps}}^{2}}{2\eps}\Big)
\end{split}
\]
has a maximum at $(\hat x, \hat t)$. Thus, in particular
\begin{equation}
\label{eq:time-deriv}
\begin{split}
\partial_t u_{\eps}(\hat x, \hat t)=\partial_t \vp (\hat x, \hat t)=\frac{\hat t-t_{\eps}}{\eps},\quad
D u_{\eps}(\hat x, \hat t)=D\vp(\hat x,\hat t)=(\hat x-x_{\eps})\frac{\abs{\hat x-x_{\eps}}^{{\hat q}-2}}{\eps^{{\hat q}-1}}.
\end{split}
\end{equation}
Moreover, since $Du_{\eps}(\hat x, \hat t)=0$, it follows that $x_{\eps}=\hat x$ and  by the definition of the inf-convolution
\[
\begin{split}
u(x,t)+\frac{\abs{\hat x-x}^{\hat q}}{{\hat q}\eps^{{\hat q}-1}}+\frac{\abs{\hat t-t}^{2}}{2\eps}\ge u_{\eps}(\hat x, \hat t)=u(\hat x,t_{\eps})+\frac{\abs{\hat t-t_{\eps}}^{2}}{2\eps}.
\end{split}
\]
Arranging the terms as
\[
\begin{split}
u(x,t)\ge u(\hat x,t_{\eps})+\frac{\abs{\hat t-t_{\eps}}^{2}}{2\eps}-\frac{\abs{\hat x-x}^{\hat q}}{{\hat q}\eps^{{\hat q}-1}}-\frac{\abs{\hat t-t}^{2}}{2\eps}
\end{split}
\]
we see that
\[
\begin{split}
\phi(x,t):= u(\hat x,t_{\eps})+\frac{\abs{\hat t-t_{\eps}}^{2}}{2\eps}-\frac{\abs{\hat x-x}^{\hat q}}{{\hat q}\eps^{{\hat q}-1}}-\frac{\abs{\hat t-t}^{2}}{2\eps}
\end{split}
\]
touches $u$ at $(\hat x,t_{\eps})$ from below. Moreover, since $u$ is a viscosity solution and
$\lim_{\hat x\neq x\to\hat x }F(D\phi(x,t_{\eps}),D^2\phi(x,t_{\eps}))=0$ since ${\hat q}> q/(q-1)$, it follows by this, $x_{\eps}=\hat x$, and  (\ref{eq:time-deriv}) that
\[
\begin{split}
0\le \partial_t \phi (\hat x, \hat t)=\frac{\hat t-t_{\eps}}{\eps}=\partial_t u_{\eps} (\hat x, \hat t)
\end{split}
\]
as claimed.

Proof of $(i)$. There exists $\vp\in C^{\infty}_0( Q_T)$ touching $u_{\eps}$ at $(\hat x,\hat t)$ from below such that $\partial_t \vp(\hat x,\hat t)=\partial_t u_{\eps}(\hat x, \hat t), D\vp(\hat x,\hat t)=Du_{\eps}(\hat x,\hat t),\ D^2\vp(\hat x,\hat t)\le D^2u_{\eps}(\hat x,\hat t)$. Then recall that the argument leading to (\ref{eq:time-deriv}) implies
\[
\begin{split}
D u_{\eps}(\hat x, \hat t)=D\vp(\hat x,\hat t)=(\hat x-x_{\eps})\frac{\abs{\hat x-x_{\eps}}^{{\hat q}-2}}{\eps^{{\hat q}-1}}
\end{split}
\]
and also
\[
\begin{split}
D^2\vp(\hat x,\hat t)\le \frac{\abs{\hat x-x_{\eps}}^{{\hat q}-2}}{\eps^{{\hat q}-1}}\Big(({\hat q}-2)\frac{\hat x-x_{\eps}}{\abs{\hat x-x_{\eps}}}\otimes \frac{\hat x-x_{\eps}}{\abs{\hat x-x_{\eps}}}+I\Big).
\end{split}
\]
The first equality above implies that if $Du_{\eps}(\hat x, \hat t)=0$ then $\hat x=x_{\eps}$, and the second since $\hat q> 2$ that $D^2\vp(\hat x,\hat t)\le 0$. If $Du_{\eps}(\hat x, \hat t)\neq 0$, then the previous inequality implies
\[
\begin{split}
D^2\vp(\hat x,\hat t)\le (\hat q-1)\frac{\abs{\hat x-x_{\eps}}^{{\hat q}-2}}{\eps^{{\hat q}-1}} I.
\end{split}
\]
\qed

\noindent {\textbf{\large \sc Acknowledgments.}} The first author was partially supported by the Academy of Finland project \#260791. Second author partially supported by Spanish Project MTM2014-52240-P.  The work was partially done while both authors visited Institut Mittag-Leffler in the fall of 2016 in the program 'Interactions between Partial Differential Equations \& Functional Inequalities'.

\def\cprime{$'$} \def\cprime{$'$} \def\cprime{$'$}


\

{
\noindent Addresses:

\noindent Mikko Parviainen\\
Department of Mathematics and Statistics, \\
University of Jyv\"askyl\"a, \\
PO~Box~35, FI-40014 Jyv\"askyl\"a, Finland\\
E-mail: mikko.j.parviainen@jyu.fi.

\

\noindent Juan Luis V\'azquez\\
Departamento de Matem\'aticas, \\
University Aut\'onoma de Madrid, \\
28049 Madrid, Spain\\
E-mail: juanluis.vazquez@uam.es.
}

\end{document}